%% file: shapereg_ims.tex
\documentclass[aos,preprint,dblspace,leqno]{imsart}

\usepackage{amsthm,amsmath,natbib}
\usepackage{amssymb}
\usepackage{amsfonts}
\usepackage{psfig}

\startlocaldefs
\newcommand{\rz}{{\mathbb R}}

\newcommand{\mathbold}[1]{\mbox{\boldmath $#1$}}

\endlocaldefs

\begin{document}

\begin{frontmatter}

\title{Nonparametric Regression, Confidence Regions and Regularization}
\runtitle{Nonparametric Regression, Confidence Regions}


\author{\fnms{P.~L.} \snm{Davies}\ead[label=e1]{laurie.davies@uni-due.de}\thanksref{t1}}
\and
\author{\fnms{A.} \snm{Kovac}\ead[label=e2]{a.kovac@bristol.ac.uk}}
\and
\author{\fnms{M.} \snm{Meise}\ead[label=e3]{monika.meise@uni-due.de}\thanksref{t1}}
\address{\printead{e1}}
\address{\printead{e2}}
\address{\printead{e3}}

\affiliation{University of Duisburg-Essen and Technical University
  of Eindhoven\\ University of Bristol\\University of Duisburg-Essen} 
\thankstext{t1}{Research supported in part by Sonderforschungsbereich 475,
  University of Dortmund} 
\runauthor{P.~L.~Davies et al.}

\begin{abstract}
In this paper we offer a unified approach to the problem
  of nonparametric regression on the unit interval. It is based on a
  universal, honest and non-asymptotic confidence region ${\mathcal
    A}_n$ which is defined by a set of linear inequalities involving
  the values of the functions at the design points. Interest will
  typically centre on certain simplest functions in ${\mathcal A}_n$
  where simplicity can be defined in terms of shape (number of local
  extremes, intervals of convexity/concavity) or smoothness (bounds on
  derivatives) or a combination of both. Once some form
  of regularization has been decided upon the confidence region can be
  used to provide honest non-asymptotic confidence bounds which are
  less informative but conceptually much simpler. 
\end{abstract}



\end{frontmatter}

\setcounter{equation}{0}
\input{shapereg1}

\input{shapereg2}

\input{shapereg4}

\bibliographystyle{apalike}
\bibliography{shapereg_ims}

\newpage
\input{shaperegapp}

\end{document}

%% file: shapereg1.tex
\section{Introduction}
Non-parametric regression on the unit interval is concerned with
specifying functions ${\tilde f}_n$ which are reasonable 
representations of a data set ${\boldsymbol
  y}_n=\{(t_i,y(t_i)),i=1,\ldots,n\}.$ The design points
$t_i$ are assumed to be ordered. Here and below we use lower case letters to
denote generic data and upper case letters to denote data generated
under a specific stochastic model. The first approach to the problem
used kernel estimators with a fixed bandwidth (Watson,
1964\nocite{WATS64}) but since then many other procedures have been 
proposed. We mention splines (Green and Silverman,
1994\nocite{GRESILV94}; Wahba, 
1990\nocite{Wahba1990}), wavelets (Donoho and Johnstone, 
1994\nocite{DOJO94}), local polynomial regression (Fan and Gijbel,
1996\nocite{FANGIJ96}), kernel estimators with local bandwidths (Wand
and Jones, 1995\nocite{WANJON95}) very often
with Bayesian and non-Bayesian versions.

The models on which the methods are based are of the form 
\begin{equation} \label{genmod}
Y(t)=f(t) +\sigma(t) \varepsilon(t), \quad t \in [0,\,1]
\end{equation}
with various assumptions being made about $\sigma(t)$, the noise
$\varepsilon(t)$ as well as the design points $\{t_1,\ldots,t_n\}.$ We
  shall restrict attention to the simplest case
\begin{equation} \label{standmod}
Y(t)=f(t) +\sigma Z(t), \quad t \in [0,\,1]
\end{equation}
where $Z$ is Gaussian white noise and the $t_i$ are given by
$t_i=i/n.$ We mention that the same ideas can be used
for the more general model (\ref{genmod}) and that robust versions are
available. The central role in this paper is played by a confidence
region ${\mathcal A}_n$ which is defined below. It specifies all
functions $\tilde f_n$ for which the model (\ref{standmod}) is consistent (in a
well--defined sense) with the data ${\boldsymbol y}_n.$ By
regularizing within ${\mathcal A}_n$ we can control both the shape and
the smoothness  of a regression function and provide honest non-asymptotic
confidence bounds. 

The paper is organized as follows. In Section
\ref{secconfreg} we define the confidence region ${\mathcal A}_n$ and
show that it is honest and non-asymptotic for data generated under
(\ref{standmod}). In Section \ref{shaperegular} we consider shape
regularization, in \ref{smooSSthregular} regularization by smoothness
as well as the combination of shape and smoothness regularization. Finally in
Section \ref{honestconbounds} we show how honest and non--asymptotic
confidence bounds can be obtained both for shape and smoothness regularization.

\section{The confidence region ${\mathcal A}_n$} \label{secconfreg}    
\subsection{Non-parametric confidence regions}\label{nonparconfreg}
 Much attention has been given to confidence sets in recent
 years. These sets are often expressed as a ball centred at some
 suitable estimate (Li, 1989\nocite{LI89}; Hoffmann and Lepski,
 2002\nocite{HOFLEP02}; Baraud, 2004\nocite{BARA04}; Cai and Low,
 2006\nocite{CAILOW06}; Robins and van der Vaart,
 2006\nocite{ROBVAA06}) with particular emphasis on adaptive methods
 where the radius of the ball automatically decreases if $f$ is
 sufficiently smooth. The concept of adaptive confidence balls is not
 without conceptual difficulties as the discussion of Hoffmann and
 Lepski (2002) shows. An alternative to smoothness is the imposition
 of shape constraints such as monotonicity and convexity (D\"umbgen,
 1998\nocite{DUEM98A}, 2003\nocite{DUEM03}; D\"umbgen and Spokoiny,
 2001; D\"umbgen and Johns, 2004\nocite{DUEMJO04}; D\"umbgen,
 2007\nocite{DUEM07}). Such confidence sets require only that $f$
 satisfy the  shape constraint which often has some independent justification.

We consider data ${\mathbold Y}_n={\mathbold Y}_n(f)$ generated under (\ref{standmod})
and limit attention to functions $f$ in some family ${\mathcal F}_n.$
We call a confidence set ${\mathcal C}_n({\mathbold Y}_n(f),\alpha)$ exact if
\begin{equation} 
P(f \in {\mathcal C}_n({\mathbold Y}_n(f),\alpha))= \alpha \quad \text{for all}\quad f
\in {\mathcal F}_n,
\end{equation}
honest (Li, 1989) if 
\begin{equation} 
P(f \in {\mathcal C}_n({\mathbold Y}_n(f),\alpha)) \ge \alpha \quad \text{for all}\quad f
\in {\mathcal F}_n,
\end{equation}
and asymptotically honest if
\begin{equation} 
\liminf_{n \rightarrow \infty}\,\,\inf_{f \in {\mathcal F}_n}
  \,\,P(f \in {\mathcal C}_n({\mathbold Y}_n(f),\alpha)) \ge \alpha
\end{equation}
holds but it is not possible to specify the $n_0$ for which the
coverage probability exceeds $\alpha -\epsilon$ for all $n \ge n_0.$
Finally we call ${\mathcal C}_n({\mathbold Y}_n(f),\alpha)$ universal
if ${\mathcal F}_n=\{f: f:[0,1] \rightarrow \rz\}.$

\subsection{Definition of ${\mathcal A}_n$}
\label{sec:confregdef}
The confidence region ${\mathcal A}_n$ we use was first given in Davies
and Kovac (2001)\nocite{DAVKOV01}. It is constructed as follows.  For
any function $g:[0,\,1] \rightarrow \rz$ and any interval
$I=[t_j,t_k]$ of $[0,1]$ with $j\le k$  we write
\begin{equation} \label{sumresid}
w({\boldsymbol y}_n,g,I)=\frac{1}{\sqrt{\vert I\vert}}\sum_{t_i \in
I}(y(t_i)-g(t_i))  
\end{equation}
where $\vert I\vert$ denotes the number of points $t_i$ in $I.$ 
With this notation
\begin{equation} \label{approxreg}
{\mathcal A}_n={\mathcal A}_n({\boldsymbol y}_n, {\mathcal I}_n,
\sigma,\tau_n)=\big\{g: \max_{ I \in {\mathcal I}_n} \vert
w({\boldsymbol y}_n,g,I)\vert \le \sigma \sqrt{\tau_n \log n\,}\,\big\}
\end{equation}
where ${\mathcal I}_n$ is a family of intervals of $[0,\,1]$ and for
given $\alpha$ the value of $\tau_n=\tau_n(\alpha)$ is defined by
\begin{equation} \label{taualpha}
P\Big( \max_{I\in {\mathcal I}_n} \frac{1}{\sqrt{\vert
    I\vert}}\Big\vert\sum_{t_i \in I} Z(t_i)\Big \vert \le 
\sqrt{\tau_n \log n\,}\,\Big)=\alpha. 
\end{equation}
If the data ${\boldsymbol y}_n$ were generated under (\ref{standmod})
then (\ref{taualpha}) implies that $P ( f \in {\mathcal A}_n)=\alpha$
with no restrictions on $f$ so that ${\mathcal A}_n$ is a universal,
exact $\alpha$--confidence region. We mention that by using an
appropriate norm (Mildenberger, 2006\nocite{MILD06}) ${\mathcal  A}_n$
can also be expressed as a ball centred at the observations
${\mathbold y}_n.$

A function $g$ belongs to ${\mathcal A}_n$ if and only if its vector
of evaluations at the design points $(g(t_1),\ldots,g(t_n))$ belongs
to the convex polyhedron in $\rz^n$ which is defined by the linear
inequalities 
\[ \frac{1}{\sqrt{\vert I\vert}}\,\Big\vert\sum_{t_i \in
I}(y(t_i)-g(t_i))\Big\vert \le \sigma_n\sqrt{\tau_n \log n}, \quad I \in
{\mathcal I}_n.\]
The remainder of the paper is in one sense nothing more than exploring the 
consequences of these inequalities for shape and smoothness
regularization. They enforce both local and global adaptivity to the
data and they are tight in that they yield optimal rates of
convergence for both shape and smoothness constraints.

 In the theoretical part of the paper we take ${\mathcal I}_n$ to be
 the set of all intervals of the form $[t_i,t_j].$ For this choice of
 ${\mathcal A}_n$ checking whether $g\in {\mathcal A}_n$ for a given $g$
 involves about $n^2/2$ linear inequalities. Surprisingly there 
exist algorithms which allow this to be done with algorithmic
complexity $\text{O}(n\log n)$ (Bernholt and Hofmeister,
2006\nocite{BERHOF06}). In practice we restrict ${\mathcal
  I}_n$ to a multiresolution scheme as follows. For some $\lambda >1$ we set 
\begin{eqnarray} 
\lefteqn{\hspace*{-1.5cm}{\mathcal
    I}_n=\left\{[t_{l(j,k)},\,t_{u(j,k)}]\,: \,l(j,k)=\lfloor 
(j-1)\lambda^k+1\rfloor,\right.} \nonumber\\
&&u(j,k)=\min\{\lfloor j\lambda^k\rfloor,n\},j=1,\ldots,\lceil
n\lambda^{-k}\rceil, \label{multischeme}\\  
&&\left.  k=1,\ldots,\lceil \log n/\log \lambda\rceil\right\}.\nonumber
\end{eqnarray}
For any $\lambda >1$ we see that ${\mathcal I}_n$ now contains
$\text{O}(n)$ intervals. For $\lambda=2$ we get the wavelet
multiresolution scheme which we use throughout the paper when doing the
calculations for explicit data sets.  If ${\mathcal I}_n$ is the set
of all possible intervals it follows from a result  of D\"umbgen and
Spokoiny (2001)\nocite{DUEMSPO01} that $\lim_{n 
  \rightarrow \infty }\tau_n=2$ whatever the value of $\alpha.$  On
the other hand for any ${\mathcal I}_n$ which contains all the
degenerate intervals $[t_j,\,t_j]$ (as will always be the case) then
$\lim_{n \rightarrow \infty }\tau_n\ge 2$ whatever $\alpha.$ In the
following we simply take $\tau_n=3$ as our default 
value. This guarantees a coverage probability of at least 
$\alpha=0.95$ for all samples of size $ n \ge 500$ and it tends
rapidly to one as the sample size increases. The exact asymptotic
distribution of $\max_{ 1 \le i < j \le n} (\sum_{l=i}^j
Z_l)^2/(j-i+1)$ has recently been derived by Kabluchko (2007)\nocite{KAB07}.

As it stands the confidence region (\ref{approxreg}) cannot be used as
it requires $\sigma.$ We use the following default estimate
\begin{equation} \label{sigman}
\sigma_n=\text{median}(\vert y(t_2)-y(t_1)\vert, \ldots,\vert
y(t_n)-y(t_{n-1})\vert)/(\Phi^{-1}(0.75)\sqrt{2})
\end{equation}
where $\Phi^{-1}$ is the inverse of the standard normal distribution
function $\Phi.$ It is seen that $\sigma_n$ is a consistent estimate
of $\sigma$ for white noise data. For data generated under
(\ref{standmod}) $\sigma_n$ is 
positively biased and consequently the coverage probability will not
decrease. Simulations show that 
\begin{equation} \label{confreg}
P\big( f \in {\mathcal A}_n({\boldsymbol Y}_n, {\mathcal
  I}_n,\sigma_n,3) \big) \ge 0.95
\end{equation}
for all $n \ge 500$ and 
\begin{equation} \label{confregasymp}
\lim_{n \rightarrow \infty}\,\inf_f\,P\big( f \in {\mathcal
  A}_n({\boldsymbol Y}_n, {\mathcal I}_n,\sigma_n,3) \big) =1.
\end{equation}
In other words ${\mathcal A}_n$ is a universal, honest and
non-asymptotic confidence region for $f.$ To separate the problem of
specifying the size of the noise from the problem of investigating the
behaviour of the procedures under the model (\ref{standmod}) we shall
always put $\sigma_n=\sigma$ for theoretical results. For real data
and in all simulations however we use the $\sigma_n$ of
(\ref{sigman}).

The confidence region ${\mathcal A}_n$ can be interpreted as the
inversion of the multiscale tests that the mean of the residuals is
zero on all intervals $I \in {\mathcal I}_n.$ A similar idea is to be
found in D\"umbgen and Spokoiny (2001)\nocite{DUEMSPO01} who invert
tests to obtain confidence regions. Their tests derive from
kernel estimators with different locations and  bandwidths  where the
kernels are chosen to be optimal for certain testing problems for
given shape hypotheses. The confidence region may be expressed in
terms of linear inequalities involving the weighted residuals with the
weights determined by the kernels. The confidence region we use corresponds to
the uniform kernel on $[0,\,1].$ Because of their multiscale
character all these confidence regions allow any lack of fit to be
localized (Davies and Kovac, 2001; D\"umbgen and Spokoiny, 2001) and
under shape regularization they automatically adapt to a
certain degree of local smoothness. Universal exact confidence regions based on the
signs of the residuals $\text{sign}(y(t_i)-g(t_i))$ rather than the
residuals themselves are to be found implicitly in Davies (1995)\nocite{DAV95} and
explicitly in D\"umbgen (2003, 2007) and D\"umbgen and Johns
(2004). These require only that under the model the errors $\varepsilon(t)$
be independently distributed  with median zero. As a consequence they do
not require an auxiliary estimate of scale such as
(\ref{sigman}). Estimates and confidence bounds based on such
confidence regions are less sensitive but much more robust.


%% file: shapereg2.tex
\section{Shape regularization and local adaptivity}  \label{shaperegular}  
\subsection{Generalities}
In this section we consider shape regularization within the confidence region
${\mathcal A}_n$. Two simple possibilities are to require that the function be
monotone or that it be convex. Although much has been written about
monotone or convex regression we are not concerned with these particular
cases. Given any data set ${\mathbold y}_n$ it is always possible to calculate
a monotone regression function, for example monotone least squares. In the
literature the assumption usually made is that the $f$ in (\ref{standmod}) is
monotone and then one examines the behaviour of a monotone 
regression function. Although this case is included in the following analysis
we are mainly concerned with determining the minimum
number of local extreme points or points of inflection required for an adequate
approximation. This is STEP 2 of Mammen (1991)\cite{MAM91B}. We shall
investigate how pronounced a peak or a point of inflection must be before it
can be detected on the basis of a sample of size $n.$ These estimates are in
general conservative but they do reflect the real finite sample
behaviour of our procedures. We shall also investigate rates of
convergence between peaks and 
points of inflection. We show that these are local in the strong sense that
the rate of convergence at a point $t$ depends only on the behaviour of $f$ in
a small neighbourhood of $t$. Furthermore we show that in a certain
sense shape regularization automatically adapts to the smoothness of $f.$  All
the calculations we perform use only the shape restrictions of the
regularization 
and the linear inequalities which determine ${\mathcal A}_n$. The mathematics
are extremely simple involving no more than a Taylor expansion and are of
no intrinsic interest. We give one such calculation in detail and refer to the
appendix for the remainder.
 
\subsection{Local extreme values} \label{locextrmval}

The simplest form of shape regularization is to minimize the number of
local extreme values subject to membership of ${\mathcal A}_n.$ We  
wish to determine this minimum number and exhibit a function in
${\mathcal A}_n$ which has this number of local extreme values. This is an 
optimization problem and the taut string algorithm of Davies and 
Kovac (2001) was explicitly developed to solve it. A short description
of the algorithm used in Kovac (2007) is given in the appendix,
section 7.3. We analyse the properties of any such solution and in
particular the ability to detect peaks or points of inflection. To do
this we consider data generated 
under the model (\ref{standmod}) and investigate how pronounced a peak of the
generating function $f$ of (\ref{standmod}) must be before it is detected on
the basis of a sample of size $n.$ We commence with the case of one local maximum
and assume that it is located at $t=1/2.$ Let $I_c$ denote an
interval which contains 1/2. For any ${\tilde f}_n$ in ${\mathcal
  A}_n$ we have 
\[ \frac{1}{\sqrt{\vert I_c\vert}} \sum_{t_i\in I_c}{\tilde f}_n(t_i) \ge
\frac{1}{\sqrt{\vert I_c\vert}} \sum_{t_i\in I_c}f(t_i)
-\sigma\sqrt{3\log n\,}+\sigma Z(I_c)\] 
and hence
\begin{equation} \label{locmax1}
\max_{t_i \in I_c} \,{\tilde f}_n(t_i)\ge \frac{1}{\vert
I_c\vert}\sum_{t_i\in I_c}f(t_i)-\sigma\frac{\sqrt{3 \log
  n}-Z(I_c)}{\sqrt{\vert I_c\vert}} 
\end{equation} 
where 
\[Z(I_c)=\frac{1}{\sqrt{\vert I_c\vert}}\sum_{t_i\in I_c}Z(t_i)
\stackrel{D}{=} N(0,1).\]
Let $I_l$ and $I_r$ be intervals to the left and right of $I_c$
respectively. A similar argument gives 
\begin{align}
\min_{t_i \in I_l} \,{\tilde f}_n(t_i)&\le \frac{1}{\vert
I_l\vert}\sum_{t_i\in I_l}f(t_i)+\sigma\frac{\sqrt{3 \log
  n}+Z(I_l)}{\sqrt{\vert I_l\vert}} \label{locmax2}\\
\intertext{and}
\min_{t_i \in I_r} \,{\tilde f}_n(t_i)&\le \frac{1}{\vert
I_r\vert}\sum_{t_i\in I_r}f(t_i)+\sigma\frac{\sqrt{3 \log
  n}+Z(I_r)}{\sqrt{\vert I_r\vert}} \label{locmax3}.
\end{align} 
If now 
\begin{eqnarray}
\lefteqn{\frac{1}{\vert I_c\vert}\sum_{t_i\in I_c}f(t_i)
  -\sigma\frac{\sqrt{3 \log n}-Z(I_c)}{\sqrt{\vert I_c\vert}}}
\hspace{2cm}\nonumber\\ 
&\ge& \max\left\{ \frac{1}{\vert
I_l\vert}\sum_{t_i\in I_l}f(t_i)+\sigma\frac{\sqrt{3 \log
  n}+Z(I_l)}{\sqrt{\vert I_l\vert}},\right.\nonumber\\
&&\hspace{2cm}\left. \frac{1}{\vert I_r\vert}\sum_{t_i\in 
I_r}f(t_i)+\sigma\frac{\sqrt{3 \log
  n}+Z(I_r)}{\sqrt{\vert I_r\vert}}
\right\}\label{locmax4} 
\end{eqnarray}
then any function in ${\mathcal A}_n$ must have a local maximum in
$I_l\cup I_c\cup  I_r.$ The random variables $Z(I_c), \,Z(I_l)$ and
$Z(I_r)$ are independently and identically distributed $N(0,1)$ random
variables. With probability at least 0.99 we have $Z(I_c) \ge -2.72,
Z(I_l) \le 2.72$ and $Z(I_r) \le 2.72$ and hence we can replace
(\ref{locmax4}) by
\begin{eqnarray}
\lefteqn{\frac{1}{\vert I_c\vert}\sum_{t_i\in I_c}f(t_i)
  -\sigma\frac{\sqrt{3 \log n}+2.72}{\sqrt{\vert I_c\vert}}}
\hspace{2cm}\nonumber\\ 
&\ge& \max\left\{ \frac{1}{\vert
I_l\vert}\sum_{t_i\in I_l}f(t_i)+\sigma\frac{\sqrt{3 \log
  n}+2.72}{\sqrt{\vert I_l\vert}},\right.\label{locmax5} \\
&&\hspace{2cm}\left. \frac{1}{\vert I_r\vert}\sum_{t_i\in 
I_r}f(t_i)+\sigma\frac{\sqrt{3 \log
  n}+2.72}{\sqrt{\vert I_r\vert}}
\right\}\nonumber
\end{eqnarray}
If we now regularize by considering those functions in ${\mathcal A}_n$ with
the minimum number of local extreme values we see that this number must be at
least one. As $f$ itself has one local extreme value and belongs to ${\mathcal
  A}_n$ with probability rapidly approaching one we see that with high
probability  the minimum number is one and that this local maximum
lies in $I_l\cup I_c\cup  I_r.$

The   condition (\ref{locmax5}) quantifies a lower bound for the power of the peak so that it
will be detected with probability of at least 0.94 on the basis of a
sample of size $n\ge 500$. The precision of the location is given by
the interval $I_l\cup I_c\cup I_r.$ We apply this to the specific function  
\begin{equation} \label{box1}
f_b(t)=b((t-1/2)/0.01)
\end{equation}
where
\begin{equation} \label{box}
b(t)=\left\{\begin{array}{ll}
1,& \,\vert t \vert \le 1\\
0,&\,\text{otherwise.}\\
\end{array}\right.
\end{equation}
We denote by $f^{*}_{bn}$ a function in ${\mathcal A}_n$ which has the
smallest number of local extreme values. As the function $f_b$ of
(\ref{box1}) lies in  ${\mathcal A}_n$ with probability rapidly
tending to one and has exactly one local extreme it follows than any
such $f^{*}_{bn}$ must have exactly one local extreme.  Suppose we wish to
detect the local maximum of $f_b$ with a precision of $\delta=0.01.$
As all points in the interval $[0.49,\,0.51]$ are in 
a sense the same local maximum  we require the local
maximum of $f_{bn}^*$ to  lie in the interval $[0.48,\,0.52].$ A short
calculation with $\sigma=1$ shows that the smallest value of $n$ for
which (\ref{locmax5}) is satisfied is approximately 19500. A small
simulation study using the taut string resulted  in the peak being
found with the prescribed accuracy in $99.6\%$ of the 10000
simulations. 

We now consider a function $f$ which has exactly one local
maximum situated in $t=1/2$ and for which 
\begin{equation} \label{boundf2}
-c_2 \le f^{(2)}(t) \le -c_1 <0, \quad t \in I_0,
\end{equation}
for some open interval $I_0$ which contains the point $t=1/2.$  We
denote by $f_n^*$ a function in ${\mathcal A}_n$ which minimizes the
number of local extremes. For large $n$ any such function $f^*_n$ will
have exactly one local extreme value which is a local maximum situated at $t^*_n$ with 
\begin{equation} \label{poslocmax}
 \vert t^*_n-1/2\vert = \text{O}_f\left(\left(\frac{\log
      n}{n}\right)^{1/5}\right).
\end{equation}
An explicit upper bound for the constant in $\text{O}_f$ in terms of
$c_1$ and $c_2$ of (\ref{boundf2}) is available. We also have
\begin{equation} \label{lbndlocmax0}
f^*_n(t^*_n) \ge f(1/2)-\text{O}_f\left(\left(\frac{\log
      n}{n}\right)^{2/5}\right)
\end{equation}
with again an explicit constant available. In the other direction
\begin{equation} \label{ubndlocmax0}
f^*_n(t^*_n) \le f(1/2)+\sigma(\sqrt{3 \log n\,}+2.4).
\end{equation}
The proofs are given in the Appendix.

More generally suppose that $f$ has a
continuous second derivative and  $\kappa$ local extreme values  situated at
$0 < t^e_1 < \ldots < t^e_{\kappa}<1$ with $f^{(2)}(t^e_k) \ne 0,
k=1,\ldots,\kappa.$ If 
$f_n^* \in {\mathcal A}_n$ now denotes a function which has the
smallest number of local extreme 
values of all functions in ${\mathcal A}_n$ it follows that with
probability tending to one  $f_n^*$ will have $\kappa$ local extreme
values located at the points  $0 < t^{*e}_{n1} < \ldots <
t^{*e}_{n\kappa}<1$ with
\begin{equation} \label{locextprec}
\vert t^{*e}_{nk}-t^e_k\vert = \text{O}_f\left(\left(\frac{\log
      n}{n}\right)^{1/5}\right), k=1,\ldots, \kappa.
\end{equation}
Furthermore if $t^e_k$ is the position of a local maximum of $f$ then
\begin{equation} \label{bndlocmax}
f^*_n(t^{*e}_{nk}) \ge f(t_k^e)-\text{O}_f\left(\left(\frac{\log
      n}{n}\right)^{2/5}\right)
\end{equation}
whereas if $t^e_k$ is the position of a local minimum of $f$ then
\begin{equation} \label{bndlocmin}
f^*_n(t^{*e}_{nk}) \le f(t_k^e)+\text{O}_f\left(\left(\frac{\log
      n}{n}\right)^{2/5}\right).
\end{equation} 
In the other direction we have 
\begin{eqnarray} 
f^*_n(t^{*e}_{nk}) &\le& f(t_k^e)+\sigma\big(\sqrt{3\log
  n}+\sqrt{3\log(8+\kappa)}\,\,\big)
\label{bndlocmaxu}\\
f^*_n(t^{*e}_{nk}) &\ge& f(t_k^e)-\sigma\big(\sqrt{3\log
  n}+\sqrt{3\log(8+\kappa)}\,\,\big). \label{bndlocminl}
\end{eqnarray}
More precise bounds cannot be attained on the basis of
monotonicity arguments alone.

\subsection{Between the local extremes}
We investigate the behaviour of $f_n^*$ between the local extremes
where $f_n^*$ is monotone. For any function $g:[0,\,1] \rightarrow \rz$ we define
\begin{equation}
\Vert g\Vert_{I,\infty} = \sup\{\vert g(t)\vert: t \in I\}.
\end{equation}
Consider a point  $t=i/n$ between two local extreme values of $f$ and
write $I_{nk}^r=[i/n,(i+k)/n]$ with $k>0.$ Then  
\begin{equation} \label{rightadapt}
f_n^*(i/n)-f(i/n)\le \min_{1\le k\le k^{*r}_n}\,
\left\{\frac{k}{n}\Vert f^{(1)}\Vert_{I_{nk}^r,\infty}+
  2\sigma\sqrt{\frac{3\log n\,}{k}}\right\} 
\end{equation}
where $k^{*r}_n$ denotes the largest value of $k$ for which $f^*_n$ is
non-decreasing on $I_{nk}^r.$  It follows from (\ref{rightadapt}) and
the corresponding inequality on the 
left that as long as $f_n^*$ has the correct global monotonicity behaviour its
behaviour at a point $t$ with $f^{(1)}(t) \ne 0$ depends only on the behaviour
of $f$ in a small neighbourhood of $t.$ In particular we have asymptotically
\begin{equation} \label{rate1}
\vert f(t)-f_n^*(t)\vert \le 3^{4/3}\sigma^{2/3}\vert
f^{(1)}(t)\vert^{1/3}\left(\frac{\log n}{n}\right)^{1/3}.
\end{equation}
Furthermore if $f^{(1)}(t)=0$ on a non-degenerate
interval $I=[t_l,t_r]$ between two local extremes then for $t_l < t <t_r$ we have
$I_l^*=[t_l,t]$ and $I^*_r=[t,t_r]$ which results in
\begin{equation} \label{rate2}
\vert f(t)-f_n^*(t)\vert \le
\frac{3^{1/2}\sigma}{\min\{\sqrt{t-t_l},\sqrt{t_r-t}\}}
\left(\frac{\log n}{n}\right)^{1/2}.
\end{equation} 
The same argument shows that if 
\[\vert f(t)-f(s)\vert \le L\vert t -s\vert^{\beta}\]
with $0< \beta \le 1$ then 
\begin{equation} \label{duemspok}
\vert f(t)-f_n^*(t)\vert \le cL^{1/(2\beta +1)}
(\sigma/\beta)^{2\beta/(2\beta+1)}(\log n/n)^{\beta/(2\beta +1)} 
\end{equation} 
where 
\[c \le (2\beta +1)3^{\beta/(2\beta +1)}\left(\frac{1}{\beta
    +1}\right)^{1/(2\beta +1)} \le 4.327.\]
Apart from the value of $c$ this corresponds to Theorem 2.2 of D\"umbgen and
Spokoiny (2001).

\subsection{Convexity and concavity} \label{concavconvex}
We now turn to shape regularization by concavity and convexity. We
take an $f$ which is differentiable with derivative $f^{(1)}$ which is strictly
increasing on $[0,\,1/2]$ and strictly decreasing on $[1/2,\,1].$ We
put $I_{nk}^c=[1/2-k/n,1/2+k/n]$, $I_{nk}^l=[t_l-k/n,t_l+k/n]$ with
$t_l+k/n<1/2-k/n$ and $I_{nk}^l=[t_r-k/n,t_r+k/n]$ with $t_r-k/n>
1/2+k/n$. Corresponding to  (\ref{locmax5}) if $f$ satisfies 
\begin{eqnarray}
\lefteqn{\min_{t \in I_{nk}^l}
f^{(1)}(t)/n-\big(2\sigma(\sqrt{3\log
  n}+2.72/\sqrt{2})\big)/k^{3/2}}\hspace{1cm}\nonumber\\
&\ge& \max \left\{\max_{t \in I_{nk}^l}
f^{(1)}(t)/n+\big(2\sigma(\sqrt{3\log n}+
2.72/\sqrt{2})\big)/k^{3/2},\right.\nonumber\\
&&\left.\max_{t \in I_{nk}^r}
f^{(1)}(t)/n+\big(2\sigma(\sqrt{3\log n}+
2.72/\sqrt{2})\big)/k^{3/2}\right\}. \label{lcrconv}
\end{eqnarray}
then it follows that with probability tending to at least 0.99 the first
derivative of every differentiable function ${\tilde f}_n \in
{\mathcal A}_n$ has at least one local maximum. Let $f_n^*$ be a
differentiable function in ${\mathcal A}_n$ whose first derivative has
the smallest number of local extreme values. Then as $f$ belongs to
${\mathcal A}_n$ with probability tending to one it follows that
$f_n^{*(1)}$ has exactly one local maximum with probability tending to
at least 0.99. Suppose now that $f$ has a continuous third derivative
and $\kappa$ points of inflection located at $0<t^i_1<\ldots
<t^i_{\kappa}$ with 
\[f^{(2)}(t^i_j)=0\, \text{ and }\,f^{(3)}(t_j^i)\ne 0,\,
j=1,\ldots,\kappa.\]
If $f^*_n$ has the smallest number of points of inflection in
${\mathcal A}_n$ then if $f \in {\mathcal A}_n$ with probability
tending to one it follows that with probability tending to one
$f_n^*$ will have $\kappa$ points of inflection located at
$0<t_{n1}^{*i}<\ldots<t_{n\kappa}^{*i}<1.$  Furthermore corresponding to
(\ref{locextprec}) we have
\begin{equation} \label{locinflprec}
\vert t^{*i}_{nk}-t^i_k\vert = \text{O}_f\left(\left(\frac{\log
      n}{n}\right)^{1/7}\right),\, k=1,\ldots, \kappa.
\end{equation}
Similarly if $t^i_k$ is a local maximum of $f^{(1)}$ then
corresponding to (\ref{bndlocmax}) we have
\begin{equation} \label{bndlocmaxf1}
f^{*(1)}_n(t^{*e}_{nk}) \ge f^{(1)}(t_k^e)-\text{O}_f\left(\left(\frac{\log
      n}{n}\right)^{2/7}\right)
\end{equation}
and if $t^i_k$ is a local minimum of $f^{(1)}$ then
corresponding to (\ref{bndlocmin}) we have
\begin{equation} \label{bndlocminf1}
f^{*(1)}_n(t^{*e}_{nk}) \le f^{(1)}(t_k^e)+\text{O}_f\left(\left(\frac{\log
      n}{n}\right)^{2/7}\right).
\end{equation} 
\subsection{Between points of inflection}
Finally  we consider the behaviour of $f_n^*$ between
the points of inflection where it is then either concave or convex. We
consider a point $t=i/n$ and suppose that $f_n^*$ is convex on
$I^r_{nk}=[i/n,(i+2k)/n].$  Corresponding to (\ref{rightadapt}) we have
\begin{equation} \label{rightadaptconv}
f_n^{*(1)}(i/n)-f^{(1)}(i/n)\le \min_{1\le k\le k^{*r}_n}\,
\left\{\frac{k}{n}\Vert f^{(2)}\Vert_{I_{nk}^r,\infty}+
  4\sigma n\sqrt{\frac{3\log n\,}{k^3}}\right\} 
\end{equation}
where $k_n^{*r}$ is the largest value of $k$ such that $f_n^*$ is
convex on $[i/n,(i+2k)/n].$ Similarly corresponding to  (\ref{leftadapt})
we have
\begin{equation} \label{rightadaptconc}
f^{(1)}(i/n) -f_n^{*(1)}(i/n)\le \min_{1\le k\le k^{*l}_n}\,
\left\{\frac{k}{n}\Vert f^{(2)}\Vert_{I_{nk}^l,\infty}+
  4\sigma n\sqrt{\frac{3\log n\,}{k^3}}\right\} 
\end{equation}
where $I_{nk}^l=[i/n-2k/n,i/n]$ and $k^{*l}_n$ is the largest value of
$k$ for which $f_n^*$ is convex on $I_{nk}^l.$ If $f^{(2)}(t) \ne 0$
we have corresponding to (\ref{rate1})
\begin{equation} \label{conrate1}
\vert f_n^{*(1)}(t)-f^{(1)}(t)\vert \le 4.36\sigma^{2/5}\vert
f^{(2)}(t)\vert^{3/5}\left(\frac{\log n}{n}\right)^{1/5}. 
\end{equation}
as $n$ tends to infinity. If $f^{(2)}(t)=0$ on the non-degenerate
interval $I=[t_l,t_r]$ then for $t_l < t <t_r$ we have corresponding
to (\ref{rate2})
\begin{equation} \label{conrate2}
\vert f_n^{*(1)}(t)-f^{(1)}(t)\vert \le
\frac{4\sqrt{3}\sigma}{\min\{(t-t_l)^{3/2},(t_r-t)^{3/2}\,\}}
\left(\frac{\log n}{n}\right)^{1/2}. 
\end{equation}
The results for $f_n^*$ itself are as follows. For a point $t$ with
$f^{(2)}(t) \ne 0$  and an interval $I_{nk}^r=[t,t+2k/n]$ where $f_n^*$ is
convex we have 
\[f_n^*(t) \le f(t)+c_1(f,t)\left(\frac{k}{n}\right)\left(\frac{\log
    n}{n}\right)^{1/5} +\frac{k^2}{2n^2}\Vert
f^{(2)}\Vert_{I_{nk}^r,\infty} +4\sigma\sqrt{\frac{3 \log n}{k}}
\]
where $c_1(f,t)=4.36\sigma^{2/5}\vert f^{(2)}(t)\vert^{3/5}.$ If we
minimize over $k$ and repeat the argument for a left interval we have
corresponding to (\ref{rate1}) 
\begin{equation} \label{conrate3}
\vert f_n^*(t)-f(t)\vert \le 11.58\sigma^{4/5}\vert
f^{(2)}(t)\vert^{1/5}\left(\frac{\log n}{n}\right)^{2/5}.
\end{equation}
Finally if $f^{(2)}(t)=0$ for $t$ in the non-degenerate interval $[t_l,t_r]$
we have corresponding to (\ref{rate2}) for $t_l < t <t_r$  
\begin{equation} \label{conrate4}
\vert f_n^*(t)-f(t)\vert \le
\frac{14\sigma}{\min\{\sqrt{t-t_l},\sqrt{t_r-t}\,\}}
\left(\frac{\log n}{n}\right)^{1/2}.
\end{equation}
If the derivative $f^{(1)}$ of $f$ satisfies $\vert f^{(1)}(t)-
f^{(1)}(s)\vert \le L\vert t-s\vert^{\beta}$ with $ 0 < \beta \le 1$ then
corresponding to (\ref{duemspok}) we have
\[\vert f_n^{*(1)}(t)-f^{(1)}(t)\vert \le cL^{3/(2\beta+3)}
(\sigma/\beta)^{2\beta/(2\beta+3)}\left(\frac{\log
    n}{n}\right)^{\beta/(2\beta+3)} \]
with 
\[ c \le
2^{\beta}\left(\frac{6\sqrt{3}}{2^{\beta}}\right)^{(\beta+2)/(2\beta+3)}+
4\sqrt{3}\beta \left(\frac{2^{\beta}}{6\sqrt{3}}\right)^{3/(2\beta+3)} \,\le \,8.78.\] 
There is of course a corresponding result for $f_n^*$ itself.



%% file: shapereg4.tex
\section{Regularization by smoothness} \label{smooSSthregular}
We turn to regularization by smoothness. 
\subsection{Minimizing total variation}
We define the total variation of the $k$th derivative of a function
$g$  evaluated at the design point $t_i=i/n$ by
\begin{equation} \label{tvk}
TV(g^k):=\sum_{i=k+2}^n \left|\Delta^{(k+1)}(g(i/n))\right|, \quad k \ge 0
\end{equation}
where
\begin{equation} \label{defdelta}
\Delta^{(k+1)}(g(i/n))=\Delta^{(1)}( \Delta^{(k)}(g(i/n)))
\end{equation}
with
\[\Delta^{(1)}(g(i/n))=n(g(i/n)-g((i-1)/n)).\]
Similarly the supremum norm $\Vert g^{(k)}\Vert_{\infty}$ is defined
by
\begin{equation} \label{supnormderiv}
\Vert g^{(k)}\Vert_{\infty}=\max_i\big\vert  \Delta^{(k)}(g(i/n))\big\vert. 
\end{equation}
Minimizing either $TV(g^k)$ or $\Vert g^{(k)}\Vert_{\infty}$ subject to
$g \in {\mathcal A}_n$ leads to a linear programming
problem. Minimizing the more traditional measure of smoothness 
\[ \int_0^1 g^{(k)}(t)^2\,dt\]
subject to $g \in {\mathcal A}_n$ leads to a quadratic programming
problem which is numerically much less stable (cf.~Davies and Meise,
2005\nocite{DAVMEI05})  so we restrict attention to minimizing
$TV(g^k)$ or $\Vert g^{(k)}\Vert_{\infty}.$ 
\begin{figure}[tb]
  \centering
  \psfig{file=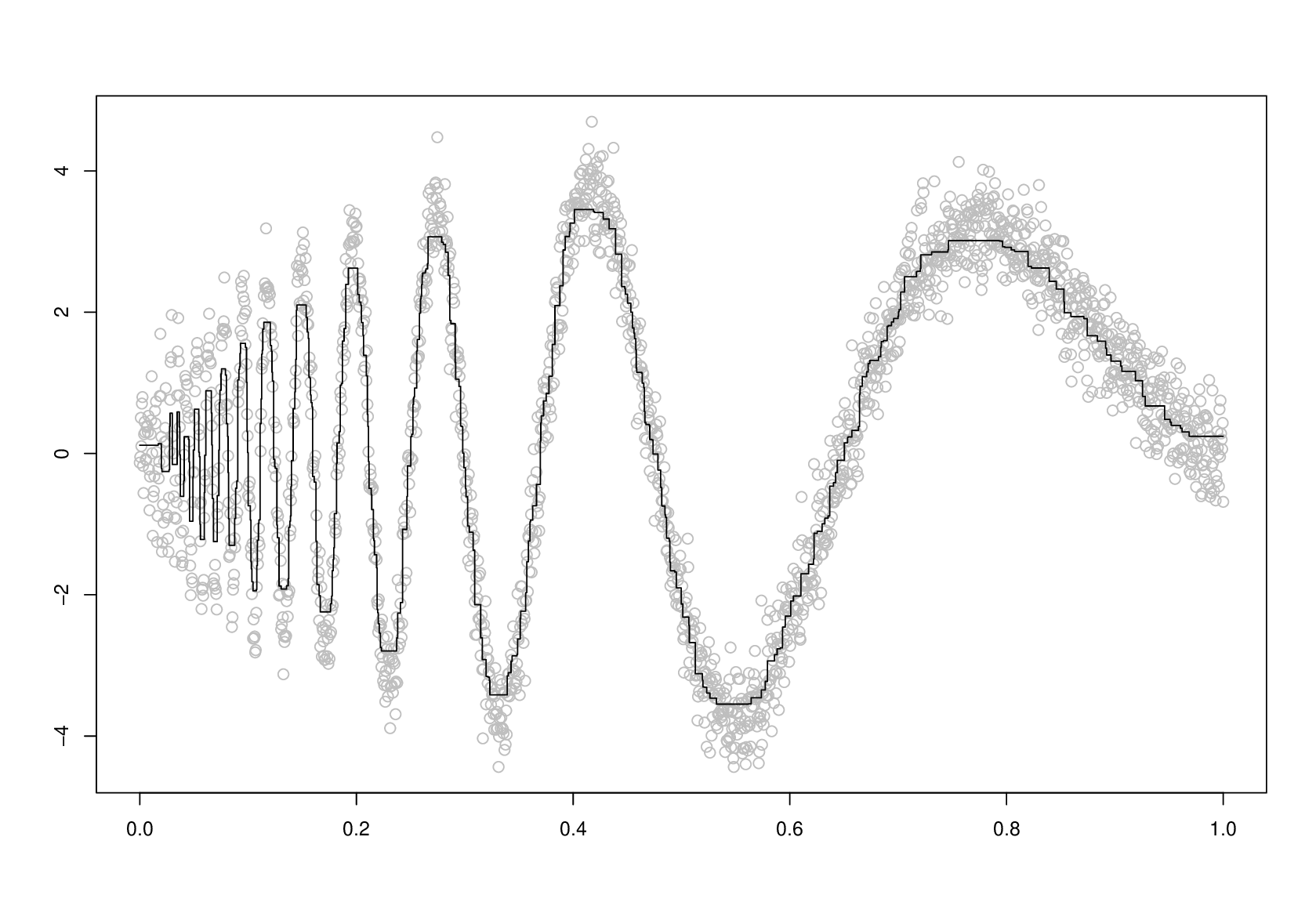,angle=270,height=6cm,width=12cm}
  \psfig{file=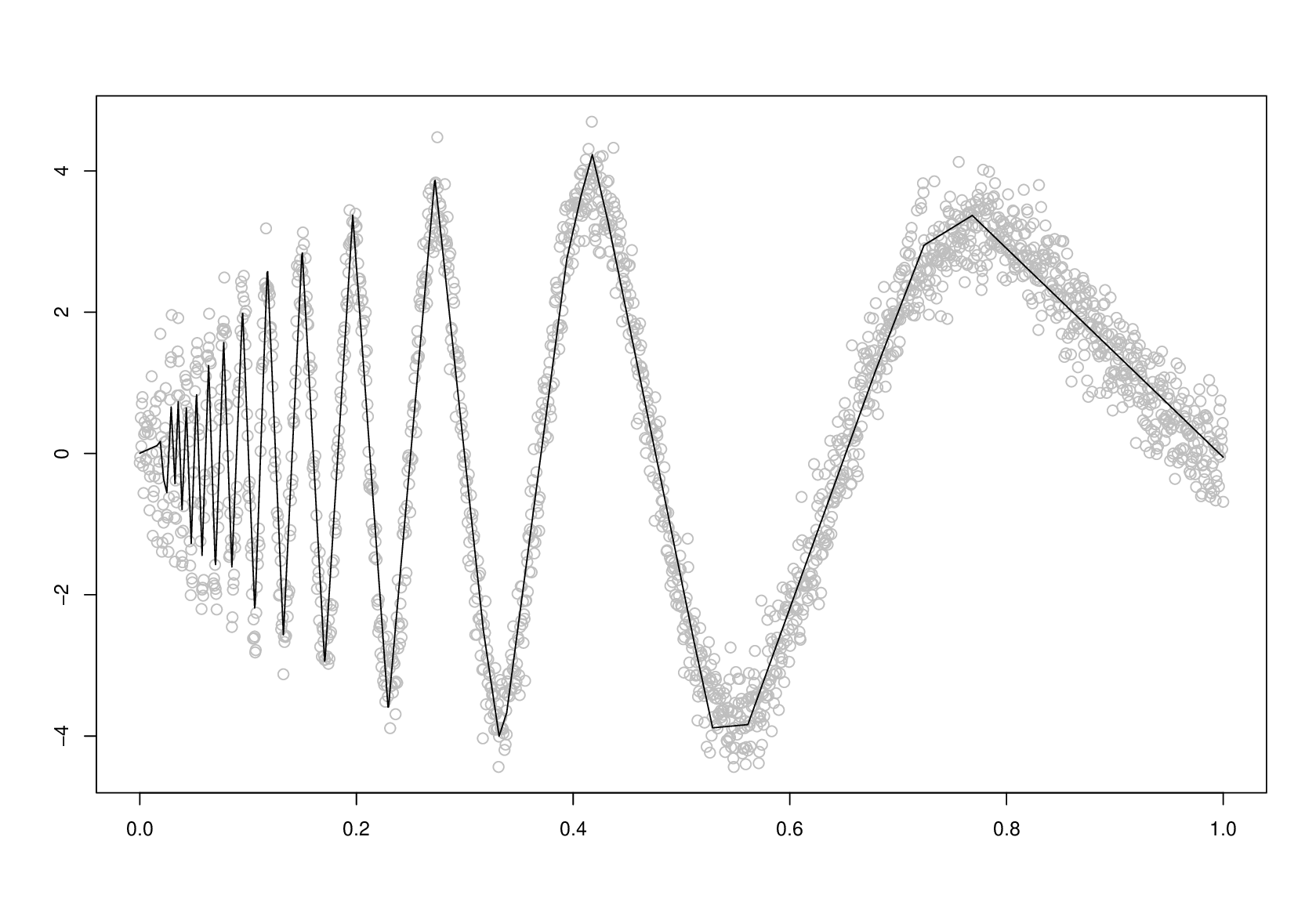,angle=270,height=6cm,width=12cm}
  \caption{ Minimization of $TV(g)$  (upper panel) and
    $TV(g^{(1)})$ (lower panel)  subject to $g \in {\mathcal
      A}_n$ for a noisy Doppler function. \label{dopplermintv}}
\end{figure}

Minimizing the total variation of $g$ itself, $k=0$, leads to piecewise
constant solutions which are very similar to the taut string
solution. In most cases the solution also minimizes the number of
local extreme values but this is not always the case. The upper
panel of Figure \ref{dopplermintv} shows the result of minimizing
$TV(g)$ for the Doppler data of Donoho and Johnstone
(1994). It has the same number of peaks as the taut string
reconstruction. The lower panel of Figure \ref{dopplermintv} shows the result
of minimizing  $TV(g^{(1)}).$  The solution is a linear
spline. Figure \ref{dopplermintv} and the following figures were
obtained using the software of Kovac (2007)\nocite{KOVFT07}. Just as
minimizing $TV(g)$ can be used for 
determining the intervals of monotonicity so we can use the solution
of minimizing $TV(g^{(1)})$ to determine the intervals of
concavity and convexity. Minimizing $TV(g^{(k)})$ or  $\Vert
g^{(k)}\Vert_{\infty}$ for larger values of $k$ leads to very smooth
functions but the numerical problems increase.

\subsection{Smoothness and shape regularization}
\begin{figure}[tb]
  \centering
  \psfig{file=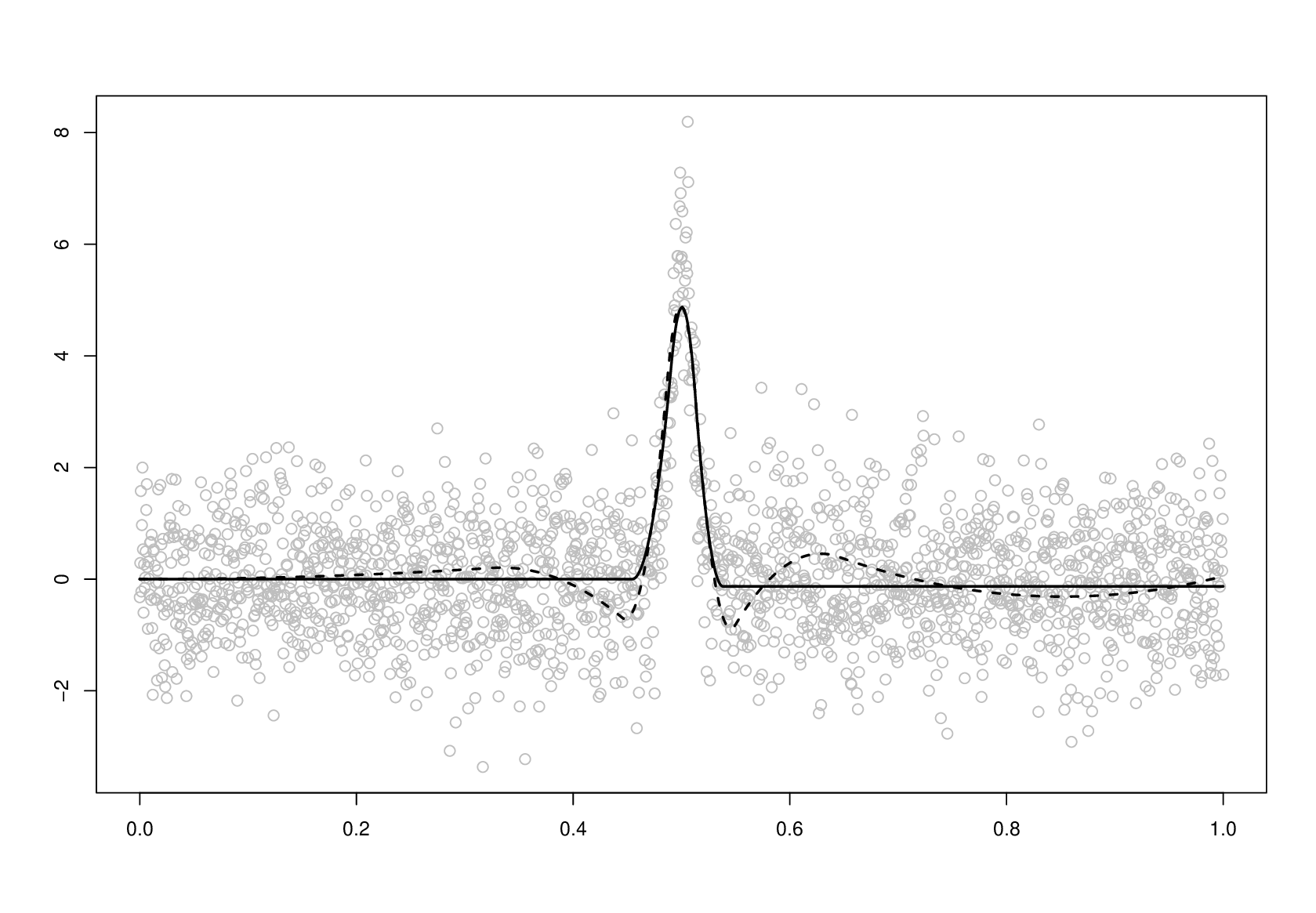,angle=270,height=6cm,width=12cm}
  \caption{ The minimization of the total variation of the second
    derivative with  (solid line) and without (dashed line) the shape 
    constraints  derived from the taut string. The solution
    subject to the shape constraints was also forced to assume the
    same value at the local maximum as the taut string solution.
 \label{mintv2mon}}
\end{figure}

Regularization by smoothness alone may lead to solutions which
do not fulfill obvious shape constraints. Figure \ref{mintv2mon} shows
the effect of minimizing the total variation of the second derivative
without further constraints and the minimization with the imposition
of the taut string shape constraints. 

\subsection{Rates of convergence} \label{smthrtconv}
 Let ${\tilde f}_n$ be such that
\begin{equation} \label{smth1}
\Vert {\tilde f}_n^{(2)}\Vert_{\infty} \le \Vert g^{(2)}\Vert_{\infty}
\quad \forall g \in {\mathcal A}_n.
\end{equation}
For data generated under (\ref{standmod}) with $f$ satisfying $\Vert
f^{(2)}\Vert_{\infty}< \infty$ it follows that with probability rapidly
tending to one 
\begin{equation} \label{smth}
\Vert {\tilde f}^{(2)}_n\Vert_{\infty} \le \Vert f^{(2)}\Vert_{\infty}.
\end{equation}
A Taylor expansion and a repetition of arguments already used leads to 
\begin{equation} \label{smthconv1}
\vert {\tilde f}_n(i/n)-f(i/n)\vert\le 3.742\Vert 
f^{(2)}\Vert_{\infty}^{1/5}\sigma^{4/5}\left(\frac{\log  
  n}{n}\right)^{2/5}
\end{equation}
on an interval 
\[\left[0.58\sigma^{2/5}(\log n)^{1/5}\big/\big(\Vert
      f^{(2)}\Vert_{\infty}^{2/5}n^{1/5}\big), \,1-0.58\sigma^{2/5}(\log
      n)^{1/5}\big/\big(\Vert
      f^{(2)}\Vert_{\infty}^{2/5}n^{1/5}\big)\right]\]
with a probability rapidly tending to one. A rate of convergence for
the first derivative may be derived in a 
similar manner and results in
\begin{equation} \label{smthconv2}
\vert {\tilde f}_n(i/n)-f^{(1)}(i/n)\vert\le 4.251\Vert 
f^{(2)}\Vert_{\infty}^{3/5}\sigma^{2/5}\left(\frac{\log  
  n}{n}\right)^{1/5}
\end{equation}
on an interval 
\[\left[2.15\sigma^{2/5}(\log n)^{1/5}\big/\big(\Vert
      f^{(2)}\Vert_{\infty}^{2/5}n^{1/5}\big), \,1-2.15\sigma^{2/5}(\log
      n)^{1/5}\big/\big(\Vert
      f^{(2)}\Vert_{\infty}^{2/5}n^{1/5}\big)\right].\]

\section{Confidence bands} \label{honestconbounds}
\subsection{The problem} \label{confboundprob}
Confidence bounds can be constructed from the confidence region
${\mathcal A}_n$ as follows. For each point $t_i$ we require a  lower bound $lb_n({\mathbold
  y}_n,t_i)=lb_n(t_i)$ and an upper bound 
$ub_n({\mathbold y}_n,t_i)=ub_n(t_i)$ such that
\begin{equation} \label{conbnds}
{\mathcal B}_n({\mathbold y}_n)=\{g: lb_n({\mathbold y}_n,t_i) \le
g(t_i) \le ub_n({\mathbold y}_n,t_i), i =1,\ldots n\}
\end{equation}
is an honest non-asymptotic confidence region
\begin{equation} \label{honconbnds} 
P( f \in {\mathcal B}_n({\mathbold Y}_n(f)))  \ge \alpha \quad \text{for
  all}\quad f \in {\mathcal F}_n
\end{equation}
for data ${\mathbold Y}_n(f)$ generated under (\ref{standmod}). In a
sense the problem has a simple solution. If we put 
\begin{equation} \label{conbnds1}
lb_n(t_i)=y(t_i)-\sigma_n\sqrt{3\log n\,},\quad ub_n(t_i)=
y(t_i)+\sigma_n\sqrt{3\log n\,},
\end{equation}
then ${\mathcal A}_n \subset {\mathcal B}_n$ and (\ref{honconbnds})
for all holds with ${\mathcal F}_n=\{f\,\vert\,f:[0,1]\rightarrow
\infty\}.$ Such universal bounds are too wide to be of any practical
use and are consequently not acceptable. They can
only be made tighter by restricting ${\mathcal F}_n$ by imposing shape or 
quantitative smoothness constraints. A qualitative smoothness
assumption such as
\begin{equation} \label{bndedf2}
{\mathcal F}_n= \{f: \Vert f^{(2)}\Vert_{\infty} < \infty\}
\end{equation}
does not lead to any improvement of the bounds (\ref{conbnds1}). They can only
be improved by replacing (\ref{bndedf2}) by a quantitative assumption such as
\begin{equation} \label{bndedf3}
{\mathcal F}_n= \{f: \Vert f^{(2)}\Vert_{\infty} < 60\}.
\end{equation}

\subsection{Shape regularization}
\subsubsection{Monotonicity}
As an example of a shape restriction we consider bounds for
non-decreasing approximations. If we denote the set of non-increasing
functions on $[0,\,1]$ by 
\[{\mathcal M}^+=\{g:\, g:[0,\,1] \rightarrow \rz,\,\, g \,\text{
  non-decreasing}\}\]
then there exists a non-decreasing approximation if and only if
\begin{equation} \label{existnondecrapprox}
 {\mathcal M}^+\cap {\mathcal A}_n \ne \emptyset.
\end{equation}
This is the case when the set of linear inequalities which define
${\mathcal A}_n$ together with $g(t_1) \le \ldots \le g(t_n)$
are consistent. This is once again a linear programming problem. If
(\ref{existnondecrapprox}) holds then the lower and upper bounds are
given respectively by
\begin{eqnarray}
lb_n(t_i) =\min\,\{g(t_i): g\in {\mathcal M}^+\cap
  {\mathcal A}_n\,\},\label{nondecrlow}\\
ub_n(t_i) =\max\,\{g(t_i): g \in {\mathcal M}^+\cap
  {\mathcal A}_n\,\}.\label{nondecrupp}
\end{eqnarray}
The calculation of $lb_n(t_i)$ and $ub_n(t_i)$ requires solving a
linear programming problem and although this can be done it is
practically impossible for larger sample sizes using standard software
because of exorbitantly long calculation times. If the family of
intervals ${\mathcal I}_n$ is restricted to a wavelet multiresolution
scheme then samples of size $n=1000$ can be handled. Fast honest
bounds can be attained as follows. If $g \in {\mathcal M}^+\cap
{\mathcal A}_n$ then for any $i$ and $k$ with 
$i+k \le n$  it follows that
\[\sqrt{k+1\,}\,\,g(t_i)\ge\frac{1}{\sqrt{k+1\,}} \,\sum_{ j=0}^k
Y_n(t_{i-j})-\sigma \sqrt{3 \log n\,}.\]
From this we may deduce the lower bound
\begin{equation} \label{fastlowerincr}
lb_n(t_i) = \max_{0 \le k \le
  i-1}\left(\frac{1}{k+1\,}\sum_{j=0}^k Y_n(t_{i-j})-\sigma\sqrt{\frac{3 \log
      n\,}{k+1}\,}\,\right)
\end{equation} 
with the corresponding upper bound
\begin{equation} \label{fastupperincr}
ub_n(t_i) = \min_{0 \le k \le
  n-i}\left(\frac{1}{k+1\,}\sum_{j=0}^k Y_n(t_{i+j})+\sigma\sqrt{\frac{3 \log
      n\,}{k+1}\,}\,\right).
\end{equation} 
Both these bounds are of algorithmic complexity $\text{O}(n^2).$  Faster
bounds can be obtained by putting
\begin{eqnarray} 
lb_n(t_i)&=& \max_{0 \le \theta(k) \le
  i-1}\left(\frac{1}{\theta(k)+1\,}\sum_{j=0}^{\theta(k)}
Y_n(t_{i-j})-\sigma\sqrt{\frac{3 \log
    n\,}{\theta(k)+1}\,}\,\right)\label{superfastlowerincr}\\
ub_n(t_i)&=& \min_{0 \le \theta(k) \le
  n-i}\left(\frac{1}{\theta(k)+1\,}\sum_{j=0}^{\theta(k)}
  Y_n(t_{i+j})+\sigma\sqrt{\frac{3 \log
      n\,}{\theta(k)+1}\,}\,\right)\label{superfastupperincr}
\end{eqnarray} 
where $\theta(k) =\lfloor \theta^k-1\rfloor$ for some $\theta >1.$
These latter bounds are of algorithmic complexity $\text{O}(n\log n).$
The fast bounds are not necessarily non--decreasing but can be made so by
putting
\begin{eqnarray*}
ub_n(t_i)&=&\min\,(ub_n(t_i),ub_n(t_{i+1})),\quad i=n-1,\ldots,1,\\
lb_n(t_i)&=&\max\,(lb_n(t_i),lb_n(t_{i-1})),\quad i=2,\ldots,n.
\end{eqnarray*}

\begin{figure}[tb]
  \centering
  \psfig{file=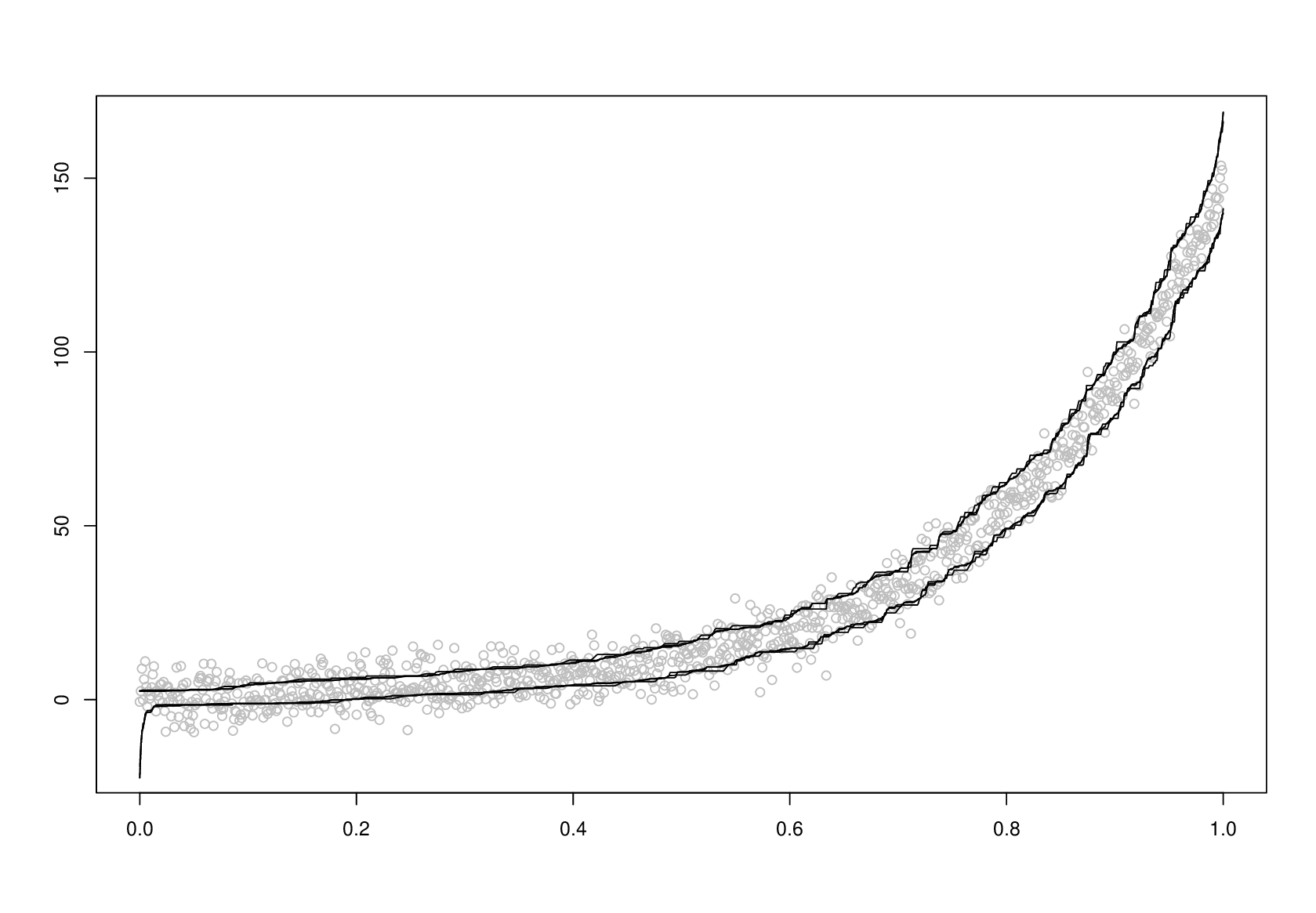,angle=270,height=5cm,width=12cm}
\psfig{file=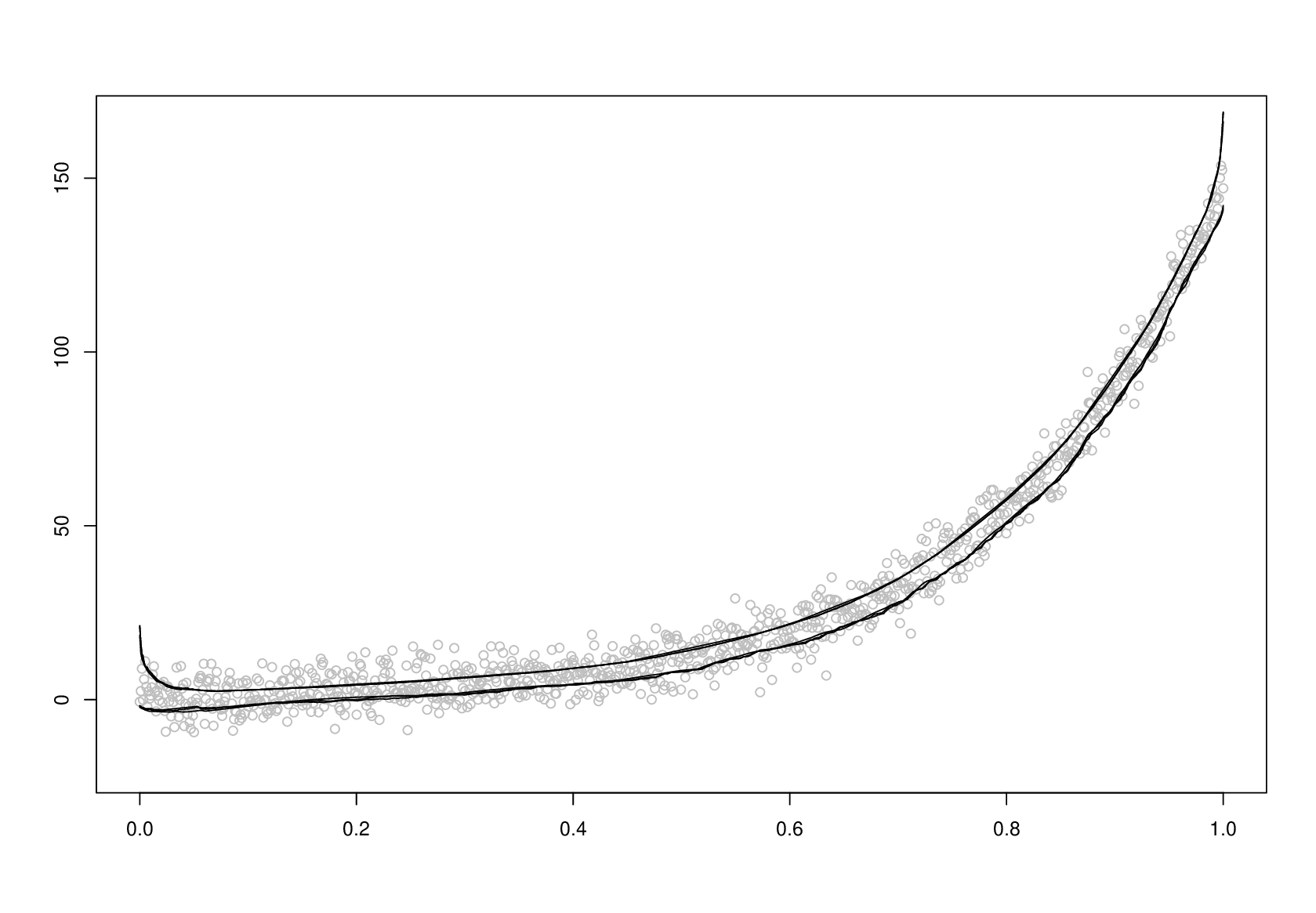,angle=270,height=5cm,width=12cm}
  \caption{ The function $f(t)=\exp(5t)$ degraded with $N(0,25)$ noise
    together with monotone confidence bounds (upper panel) and convex
    confidence bounds (lower panel). The three lower bounds in the
    upper panel are derived from (\ref{nondecrlow}), (\ref{fastlowerincr}) and
    (\ref{superfastlowerincr}) and the 
corresponding upper bounds are (\ref{nondecrupp}), (\ref{fastupperincr}) and
(\ref{superfastupperincr}). The lower bounds for the lower panel are
(\ref{lowcnvx}), (\ref{fastlowercnvx}) and (\ref{superfastlowercnvx})
and the corresponding upper bounds (\ref{uppcnvx}), (\ref{fastuppercnvx}) and
(\ref{superfastuppercnvx}) \label{expfnbndsinccnvx}}
\end{figure}

The upper panel of Figure \ref{expfnbndsinccnvx} shows data generated by 
\begin{equation} \label{expdata}
Y(t)=\exp(5t)+5Z(t)
\end{equation}
evaluated on the grid $t_i=i/1000, i=1,\ldots,100$ together with the three
lower and three upper bounds with $\sigma$ replaced by $\sigma_n$ of
(\ref{sigman}). The lower bounds are those given by (\ref{nondecrlow}) with
${\mathcal I}_n$ a dyadic multiresolution scheme,
(\ref{fastlowerincr}) and (\ref{superfastlowerincr}) with $\theta=2.$
The times required for were about 12 hours, 19 seconds and less than one second
respectively with corresponding times for the upper
bounds (\ref{nondecrupp}), (\ref{fastupperincr}) and
(\ref{superfastupperincr}). The differences between the bounds are not very
large and it is not the case that one set of bounds dominates the others. The
methods of Section \ref{shaperegular} can be applied to show that all
the uniform bounds are optimal in terms of rates of convergence.

\subsubsection{Convexity}
Convexity and concavity can be treated similarly. If we denote the
set of convex functions on $[0,\,1]$ by ${\mathcal C}^+$ then there
exists a convex approximation if and only if
\[  {\mathcal C}^+\cap {\mathcal A}_n \ne \emptyset.\]
Assuming that the design points are of the form $t_i=i/n$ this will be
the case if and only if the set of linear constraints
\[ g(t_{i+1})-g(t_i) \ge g(t_i)-g(t_{i-1}), \,\,i=2,\ldots, n-1,\]
are consistent with the linear constraints which define ${\mathcal
  A}_n.$ Again this is a linear programming problem.  If this is the
case then lower and upper bounds are given respectively by
\begin{eqnarray}
lb_n(t_i) &=&\min\,\{g(t_i): g \in {\mathcal C}^+\cap
  {\mathcal A}_n\,\},\label{lowcnvx}\\
ub_n(t_i) &=&\max\,\{g(t_i): g \in {\mathcal C}^+\cap
  {\mathcal A}_n\,\}\label{uppcnvx}
\end{eqnarray}
which again is a linear programming problem which can only be solved
for relatively small values of $n$. An honest but faster upper bound
can be obtained by noting that
\[ g(i/n) \le \frac{1}{2k+1}\sum_{j=-k}^k g((i+j)/n),\quad k \le
\min\,(i-1,n-i)\]
which gives rise to
\begin{equation} \label{fastuppercnvx}
ub_n(t_i) = \min_{0 \le k \le \min\,(i-1,n-i)}
\left(\frac{1}{2k+1\,}\sum_{j=-k}^k Y_n(t_{i+j})+\sigma\sqrt{\frac{3
      \log  n\,}{2k+1}\,}\,\right).
\end{equation} 
A fast lower bound is somewhat more complicated. Consider a function
${\tilde f}_n \in {\mathcal C}^+\cap{\mathcal A}_n$ and two points
$(i/n,{\tilde f}_n(i/n))$ and $((i+k)/n,ub_n((i+k)/n)).$ As ${\tilde
  f}_n((i+k)/n)\le ub_n((i+k)/n)$ and ${\tilde f}_n$ is convex it
follows that ${\tilde f}_n$ lies below the line joining $(i/n,{\tilde
  f}_n(i/n))$ and $((i+k)/n,ub_n((i+k)/n)).$ From this and ${\tilde
  f}_n \in {\mathcal A}_n$ we may derive a lower bound by noting
\begin{eqnarray} 
\lefteqn{lb_n(t_i)\le lb_n(t_i,k):=}\label{fastlowauxcnvx}\\
&&\max_{1 \le j \le k}\left(\frac{1}{j}\sum_{l=1}^jY_n(t_{i+j})
  -ub_n(t_{i+k})(j+1)/(2k)-\sigma\sqrt{3\log n/j}\right)\nonumber
\end{eqnarray} 
for all $k,\, -i+1\le k \le n-i.$ An honest lower bound is therefore given by
\begin{equation} \label{fastlowercnvx}
lb_n(t_i)= \max_{-i+1\le k \le n-i}\, lb_n(t_i,k).
\end{equation}
The algorithmic complexity of $ub_n$ as given by (\ref{fastuppercnvx})
is $\text{O}(n^2)$ whilst that of the lower bound (\ref{fastlowercnvx}) is
$\text{O}(n^3).$ Corresponding to (\ref{superfastupperincr}) we have
\begin{equation} \label{superfastuppercnvx}
ub_n(t_i) = \min_{0 \le \theta(k) \le \min\,(i-1,n-i)}
\left(\frac{1}{2\theta(k)+1\,}\sum_{j=-\theta(k)}^{\theta(k)}
  Y_n(t_{i+j})+\sigma\sqrt{\frac{3 \log  n\,}{2\theta(k)+1}\,}\,\right)
\end{equation} 
and to (\ref{superfastlowerincr})
\begin{equation} \label{superfastlowercnvx}
lb_n(t_i)= \max_{-i+1\le \theta(k) \le n-i}\, lb_n(t_i,\theta(k)).
\end{equation}
where
\begin{eqnarray} 
\lefteqn{\hspace{-2cm}lb_n(t_i)\le lb_n(t_i,\theta(k)):=} \hspace{-3cm}
\label{superfastlowauxcnvx}\\ 
&&\max_{1 \le \theta(j) \le
  \theta(k)}\bigg(\,\frac{1}{\theta(j)}\sum_{l=1}^{\theta(j)}Y_n(t_{i+j}) 
  -ub_n(t_{i+\theta(k)})(\theta(j)+1)/(2\theta(k))\bigg.-\nonumber\\
&&\hspace{7cm}\bigg. \sigma\sqrt{3\log n/\theta(j)\,}\,\bigg) \nonumber
\end{eqnarray} 
with $\theta(k)=\lfloor \theta^k\rfloor$ for some $\theta >1.$ The
algorithmic complexity of (\ref{superfastuppercnvx}) is $\text{O}(n\log n)$ and
that of (\ref{superfastlowercnvx}) is $\text{O}(n(\log n)^2).$

The lower panel of Figure \ref{expfnbndsinccnvx} shows the same data
as in the upper panel but with the lower bounds given by (\ref{lowcnvx}) with
${\mathcal I}_n$ a dyadic multiresolution scheme,
(\ref{fastlowercnvx}) and (\ref{superfastlowercnvx}) and the
corresponding upper bounds (\ref{uppcnvx}), (\ref{fastuppercnvx}) and
(\ref{superfastuppercnvx}). The calculation of each of the bounds
(\ref{lowcnvx}) and (\ref{uppcnvx}) took about 12 hours. The lower bound
(\ref{fastlowercnvx}) took about 210 minutes whilst
(\ref{superfastlowercnvx}) was calculated in less than 5 seconds. The
lower bound (\ref{lowcnvx}) is somewhat better  than
(\ref{fastlowercnvx}) and (\ref{superfastlowercnvx}) but the latter
two are almost indistinguishable. 

\subsubsection{Piecewise monotonicity}
\begin{figure}[h]
  \centering
 \psfig{file=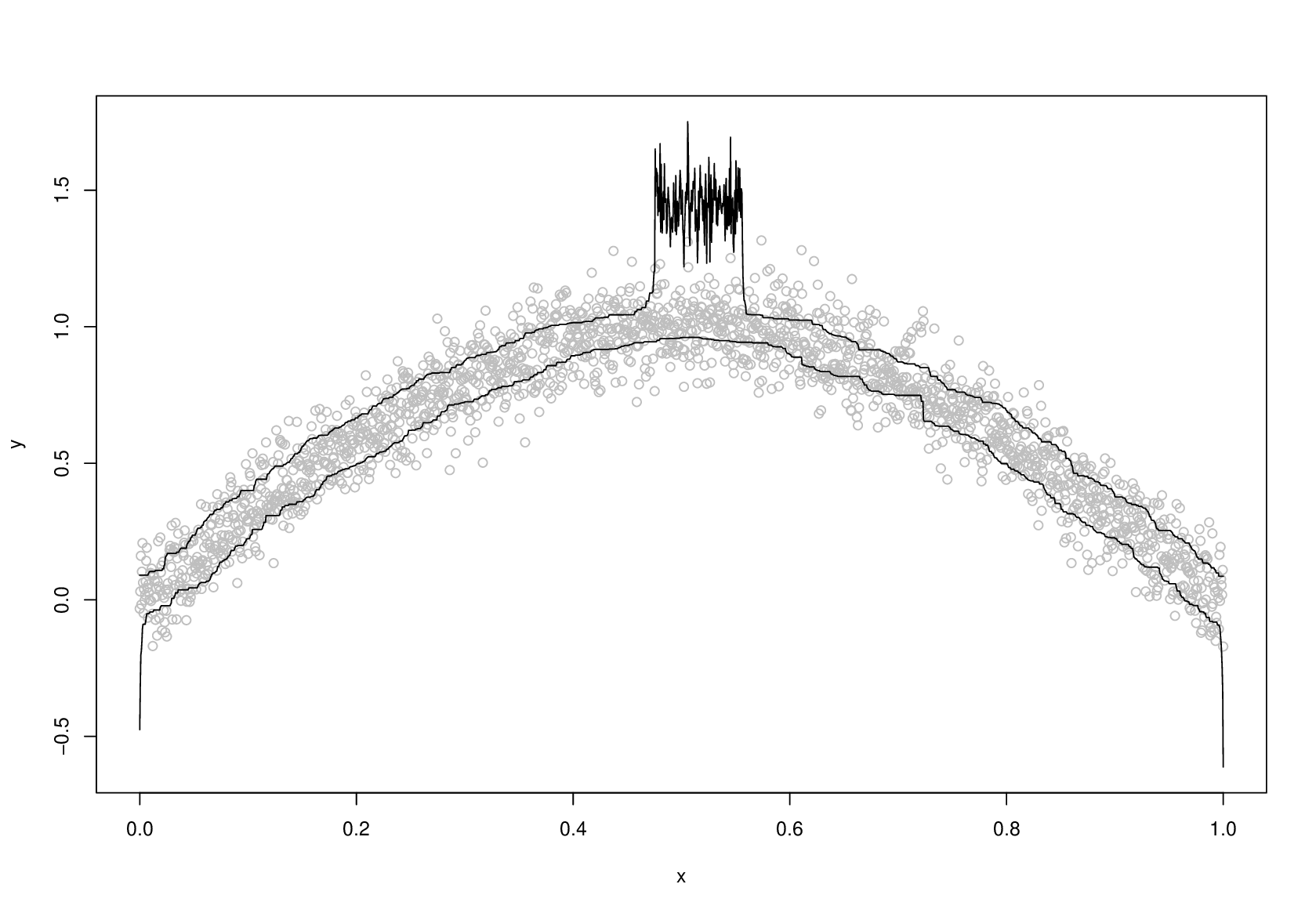,height=5cm,width=12cm,angle=270}
  \psfig{file=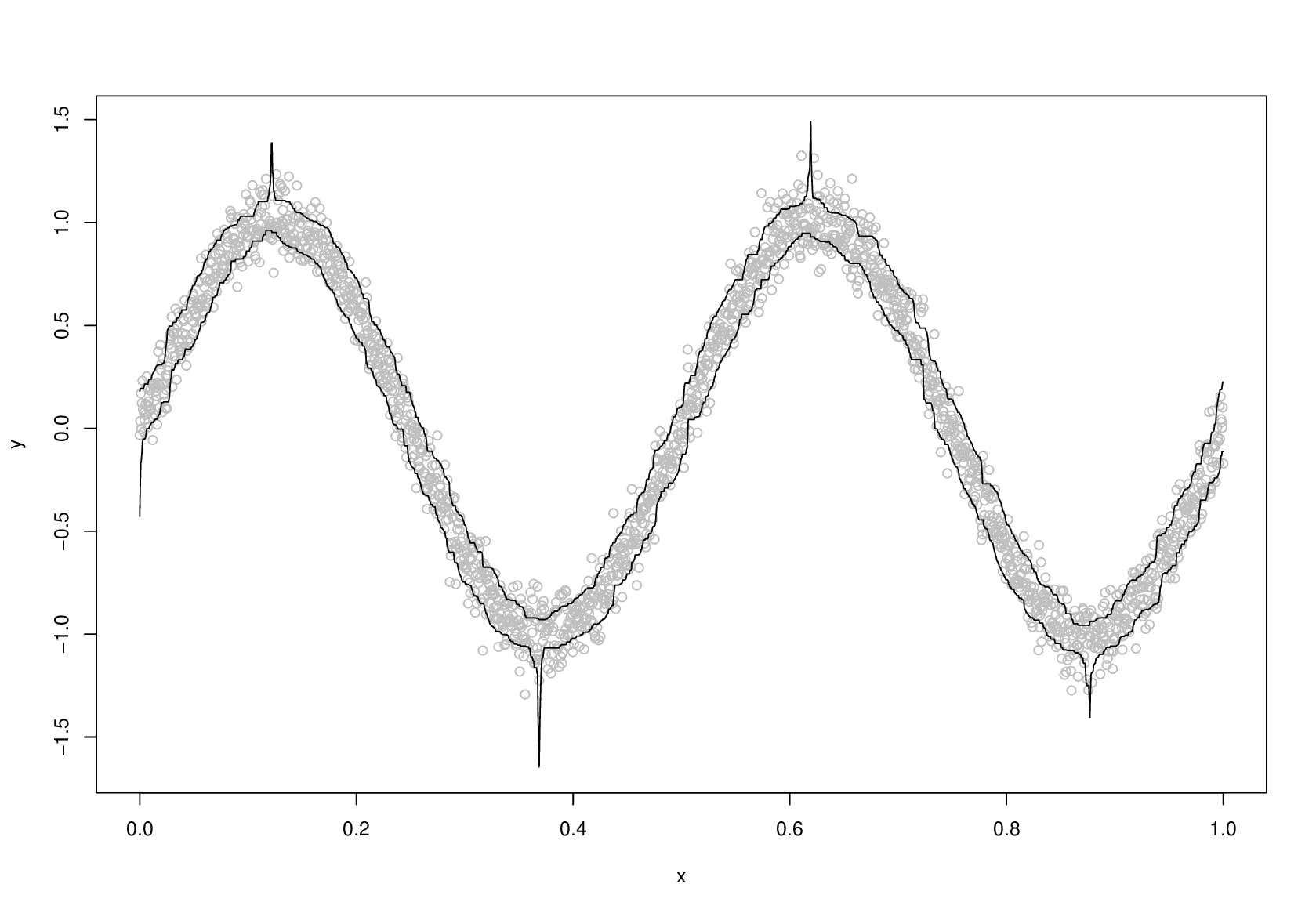,height=5cm,width=12cm,angle=270} 
  \caption{ Confidence bounds without (upper panel) and with (lower
    panel) the specification of the precise positions of the local extreme
    values. The positions in the lower panel are the default choices
    obtained from the taut string reconstruction (Kovac 2007). The
    bounds are the fast bounds (\ref{superfastlowerincr}) and
    (\ref{superfastupperincr}) with $\theta=1.5$.
    \label{EXMPLcbdav_sin_2_maxautots}}
\end{figure}

We now turn to the case of functions which are piecewise monotone. The
possible positions of the local extremes can in theory be determined
by solving the appropriate linear programming problems. The taut
string methodology is however extremely good and very fast so we can
use this solution to identify possible positions of the local
extremes. The confidence bounds depend on the exact location of the
local extreme. If we take the interval of constancy of the taut string
solution which includes the local maximum we may calculate confidence
bounds for any function which has its local maximum in this
interval. The result is shown in the top panel of Figure
\ref{EXMPLcbdav_sin_2_maxautots} where we used the fast bounds
(\ref{superfastlowerincr}) and (\ref{superfastupperincr})(\ref{superfastlowerincr}) and (\ref{superfastupperincr}) with
$\theta=1.5$. Finally if we use the mid-point of the 
taut string interval as a default choice for the position of a local
extreme we obtain confidence bounds as shown in the lower panel of Figure
\ref{EXMPLcbdav_sin_2_maxautots}. The user can of course specify these
positions and the programme will indicate if they are consistent with
the linear constraints which define the approximation region
${\mathcal A}_n.$

\begin{figure}[h]
  \centering
  \psfig{file=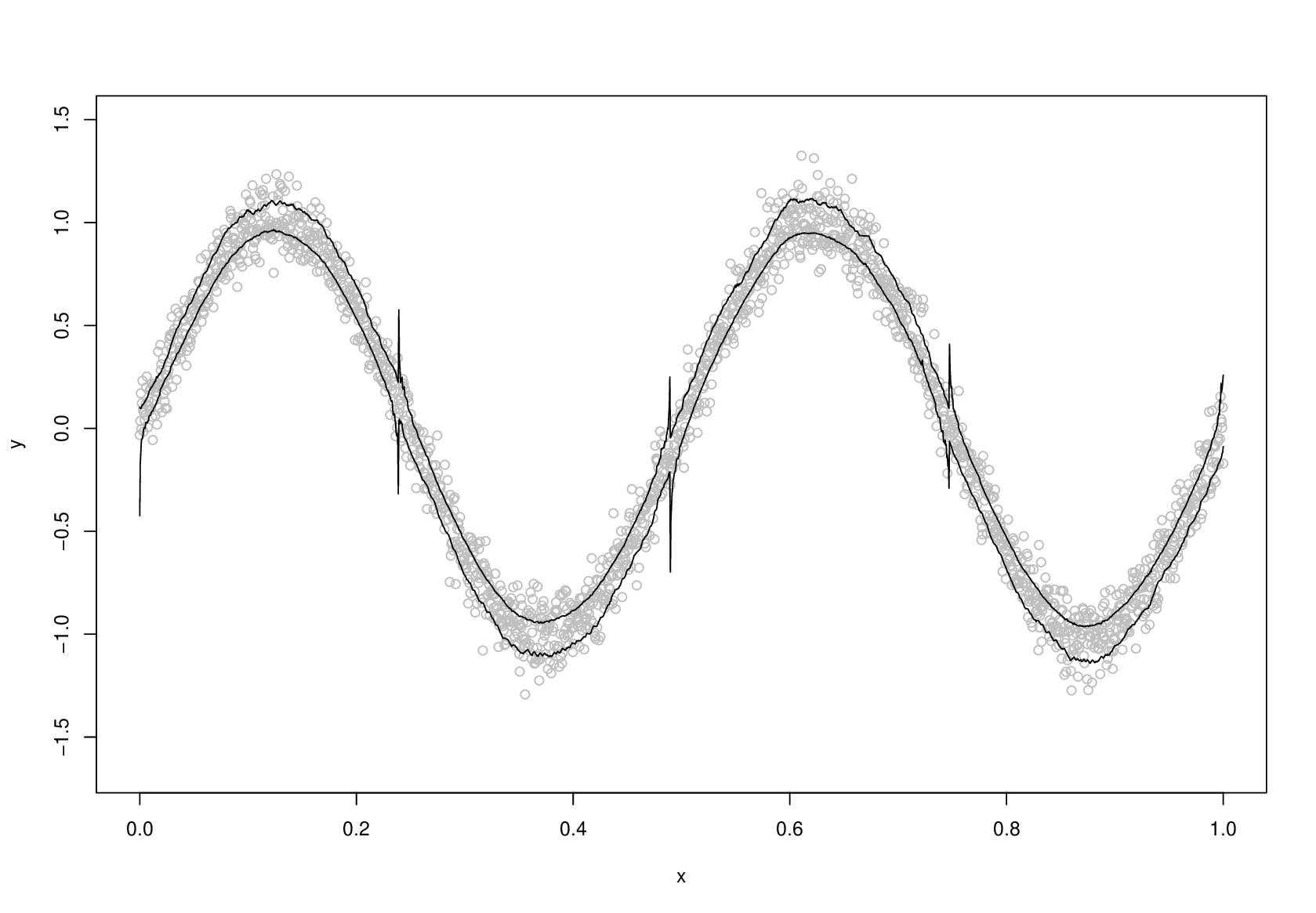,angle=270,height=5cm,width=12cm}  
  \psfig{file=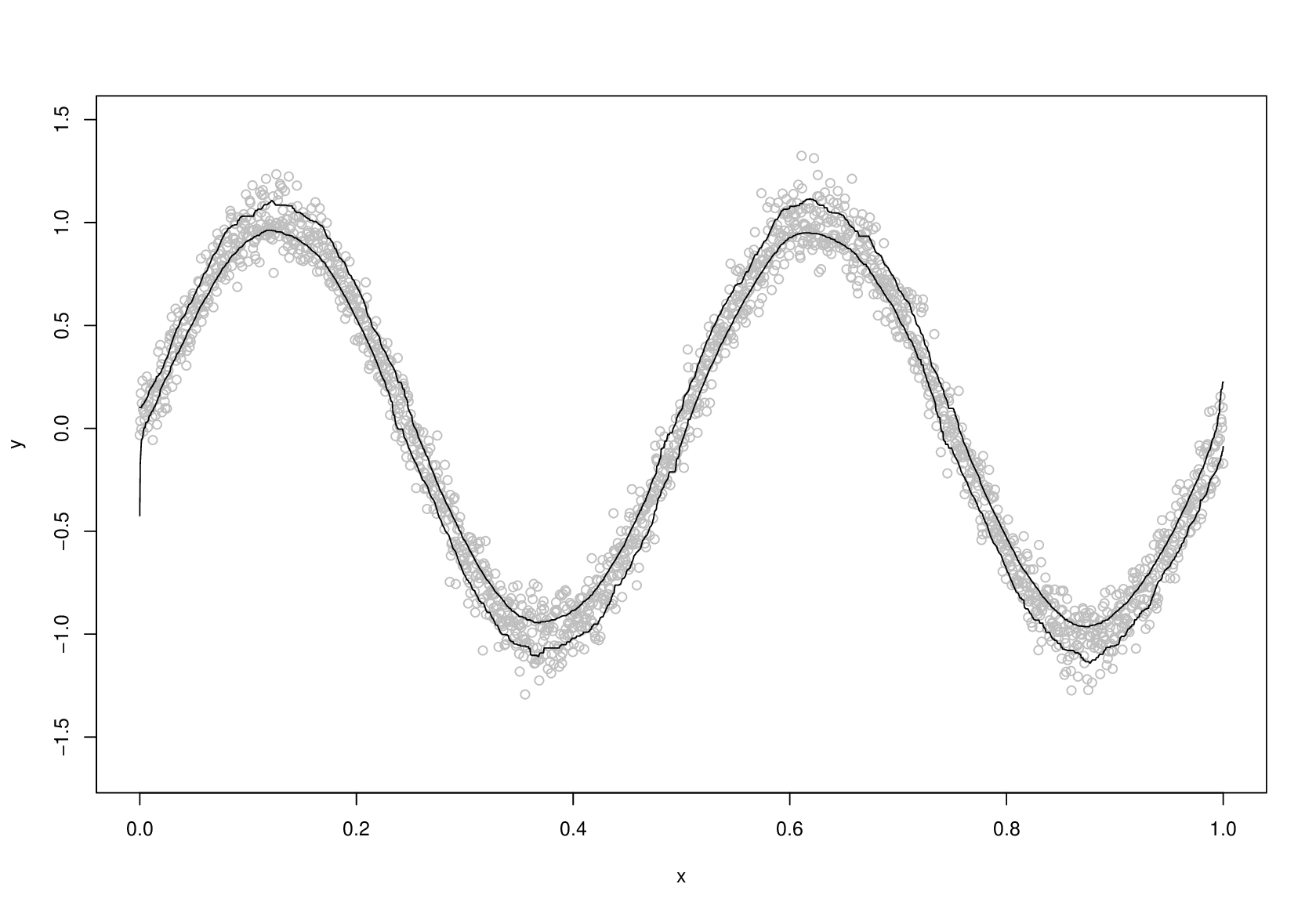,angle=270,height=5cm,width=12cm}   
  \caption{ Confidence bounds with default choices for the intervals of
    convexity/concavity (upper panel based on (\ref{superfastuppercnvx}) and
(\ref{superfastlowercnvx}) with $\theta=1.5$) and combined confidence
    bounds for default choices of intervals of monotonicity and
    convexity/concavity. \label{EXMPLcbdav_sin_2_cvx_auto}}
\end{figure}

\subsubsection{Piecewise concave--convex}
We can repeat the idea for functions which are piecewise
concave--convex. There are fast methods for determining the intervals
of convexity and concavity based on the algorithm devised by
Groeneboom (1996)\nocite{GROE96} but in this section we use the
intervals obtained by minimizing the total variation of the first
derivative  (Kovac, 2007). The upper panel of Figure
\ref{EXMPLcbdav_sin_2_cvx_auto} shows the result for
convexity/concavity which corresponds to Figure
\ref{EXMPLcbdav_sin_2_maxautots}. Finally the lower panel of Figure
\ref{EXMPLcbdav_sin_2_cvx_auto} shows the result of imposing both
monotonicity and convexity/concavity constraints. In both cases the
bounds used are the fast bounds (\ref{superfastuppercnvx}) and
(\ref{superfastlowercnvx}) with $\theta=1.5$.

\subsubsection{Sign-based confidence bounds}
As mentioned in Section \ref{sec:confregdef} work has been done on
confidence regions based on the signs of the residuals. These can also
be used to calculate confidence bands for shape-restricted functions. We
refer to  Davies (1995), D\"umbgen (2003, 2007) and D\"umbgen and Johns
(2004).

\subsection{Smoothness regularization} \label{smoothregsecconf}
\begin{figure}[p]
  \centering
  \psfig{file=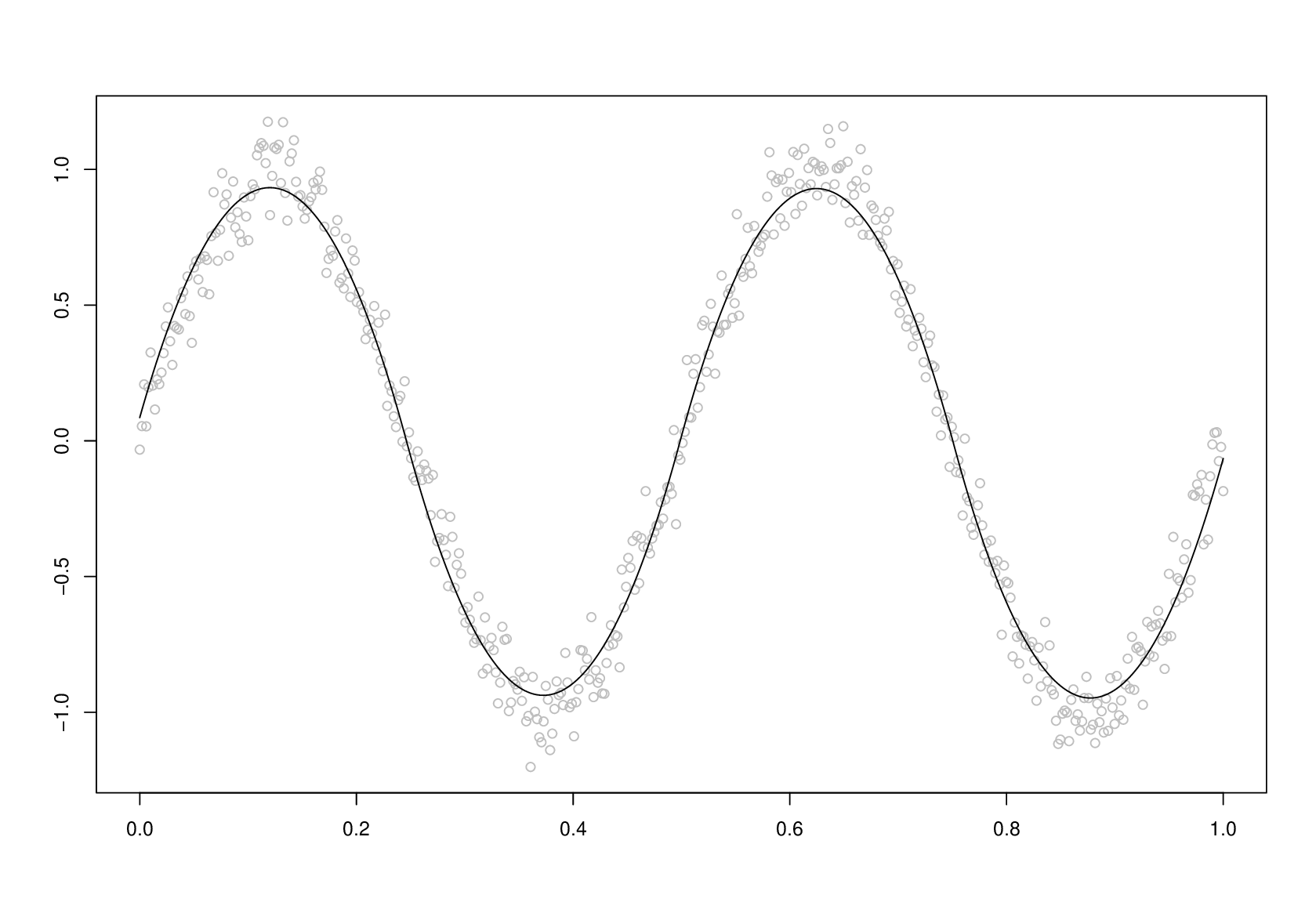,angle=270,height=5cm,width=12cm}  
  \psfig{file=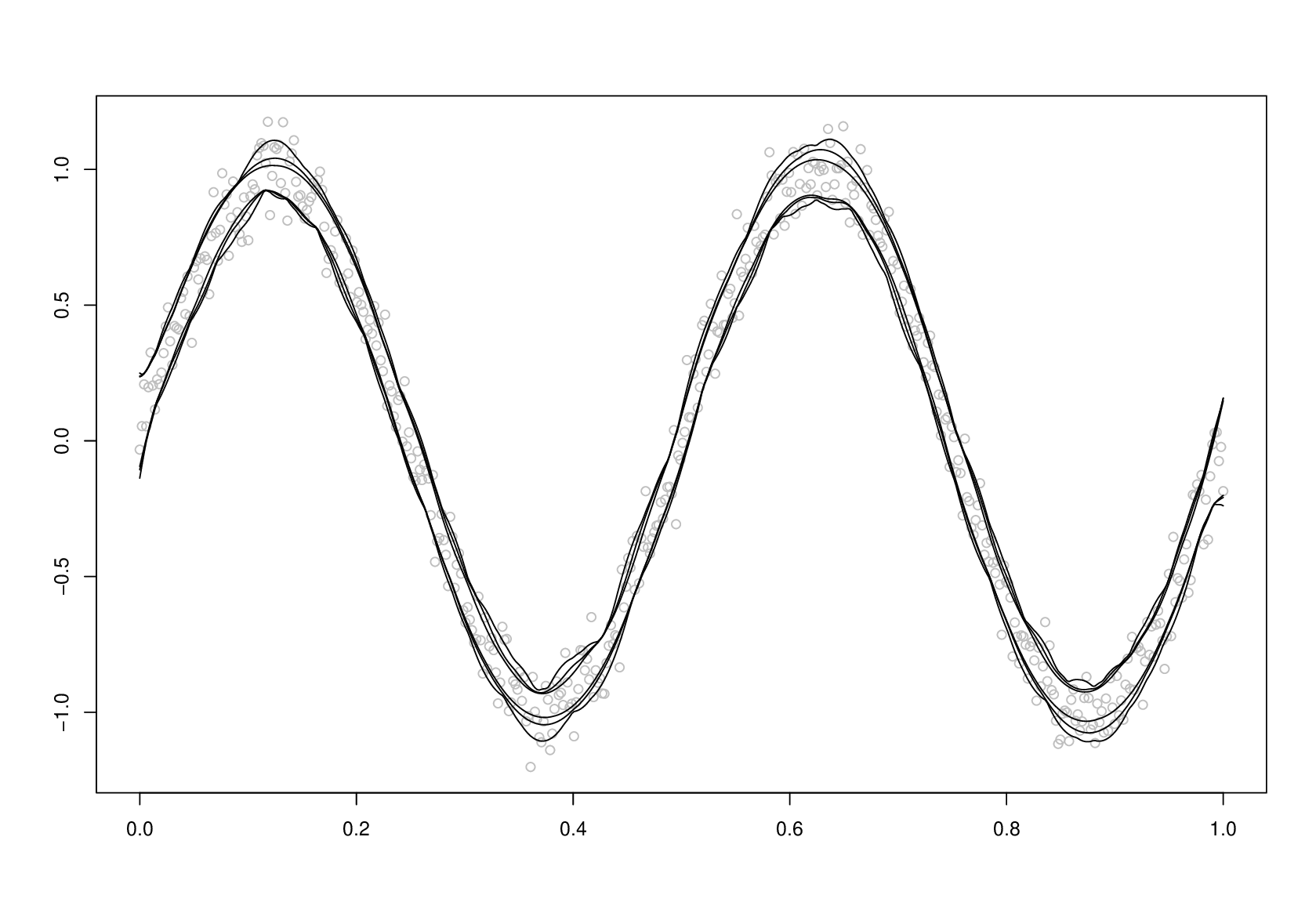,angle=270,height=5cm,width=12cm} 
  \psfig{file=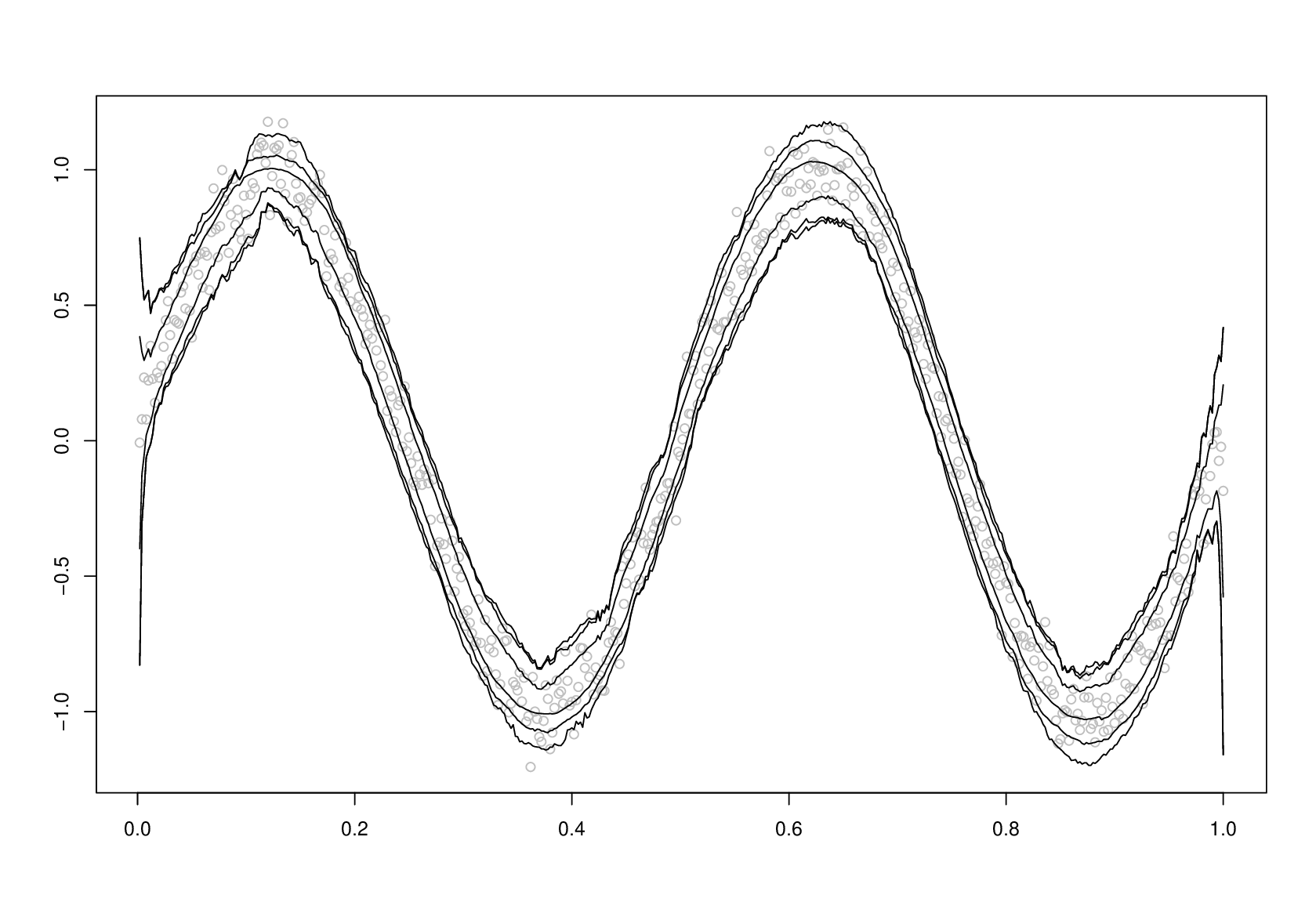,angle=270,height=5cm,width=12cm}  
  \caption{Smoothness confidence bounds for $f \in {\mathcal
      F}_n=\{f:\Vert {\tilde f}_n^{(2)}\Vert _{\infty} \le K\}$ for
    data generated according to (\ref{standmod}) with $f(t)=\sin(4\pi
    t),\,\sigma=0.2$ and $n=500.$ 
    The top panel shows the function which minimizes  $\Vert
    g^{(2)}\Vert _{\infty}$. The minimum is  117.7 compared with 
    $16\pi^2=157.9$ for $f(t)$. For this value of $K$ the bounds are
    degenerate. The centre panel shows the confidence bounds for $K=
    137.8, \,157.9$ and $315.8$  The bottom panel shows the  
  corresponding fast bounds  (\ref{lbfastsmooth}) and
  (\ref{ubfastsmooth}) with $\theta=1.5$ for the same values of $K.$ 
    \label{EXMPL_sin_cb1}}
\end{figure}

We turn to the problem of constructing lower and upper confidence
bounds under some restriction on smoothness. For simplicity we take
the supremum norm $\Vert g^{(2)}\Vert_{\infty}$ to be the measure of
smoothness for a function $g$.  The discussion in Section \ref{confboundprob}
shows that honest bounds are attainable only if we restrict $f$ to a set 
${\mathcal F}_n=\{g: \Vert g^{(2)}\Vert_{\infty} \le K\}$ with a specified
$K$. We illustrate the idea  using data
generated by (\ref{standmod}) with $f(t)=\sin(4\pi t)$ and $\sigma=1.$  The
minimum value of $\Vert g^{(2)}\Vert_{\infty}$ is 117.7 which
compares with $16\pi^2=157.9$ for $f$ itself. The upper panel of Figure
\ref{EXMPL_sin_cb1} shows  the data together with the resulting function 
$f^*_n.$  The bounds under the restriction $\Vert {\tilde
  f}_n^{(2)}\Vert_{\infty} \le 117.2$ coincide with the function $f^*_n$
itself.  The middle panel of Figure \ref{EXMPL_sin_cb1} show the bounds based
on $\Vert g^{(2)}\Vert_{\infty} \le K$ for
\[K=137.8(=(117.7+157.9)/2),\,\,157.9\quad \text{and}\quad 315.8
(=2\times 157.9).\] 
Just as before fast bounds are also available. We have for the lower
bound for given $K$
\begin{equation} \label{lbfastsmooth}
lb(i/n) \le \min_{k}\left(\frac{1}{2k+1}\sum_{j=-k}^k Y((i+j)/n)
  +\left(\frac{k}{n}\right)^2K +\sigma\sqrt{\frac{3\log
      n}{2k+1}\,}\right)
\end{equation}
and for the upper bound
\begin{equation} \label{ubfastsmooth}
ub(i/n) \ge \max_{k}\left(\frac{1}{2k+1}\sum_{j=-k}^k Y((i+j)/n)
  -\left(\frac{k}{n}\right)^2K -\sigma\sqrt{\frac{3\log
      n}{2k+1}\,}\right).
\end{equation}

As it stands the calculation of these bounds is of algorithmic
complexity $\text{O}(n^2)$ but this can be reduced to $\text{O}(n\log
n)$ by restricting $k$ to be of the form $\theta^m.$ The method also gives a
lower bound for $\Vert g^{(2)}\Vert_{\infty}$ for $g$ to be
consistent with the data. This is the smallest value of $K$ for which
the lower bound $lb$ lies beneath the upper bound $ub.$ If we do this
for the data of Figure \ref{EXMPL_sin_cb1} with $\theta=1.5$ then the smallest value is
104.5  as against the correct bound of 115.0. The lower panel of Figure
\ref{EXMPL_sin_cb1} shows the fast bounds for the same data and values of $K$.

\section{Acknowledgment}
The authors gratefully acknowledge talks with Lutz D\"umbgen which in
particular lead to the smoothness regularization described in Section
\ref{smooSSthregular}.

We also gratefully acknowledge the financial support of the German
Science Foundation (Deutsche Forschungsgemeinschaft,
Sonderforschungsbereich 475, Reduction of Complexity in Multivariate
Data Structures). We also acknowledge helpful comments made by two
referees, an Associate Editor and an Editor which have lead to a more
focussed article.


%% file: shaperegapp.tex
\section{Appendix}
\subsection{Proofs of  Section \ref{locextrmval}}
\subsubsection{Proof of (\ref{poslocmax})}
Let $k$ be such that $I_c=[1/2-k/n,\, 1/2+k/n] \subset I_0.$ A Taylor
expansion together with (\ref{boundf2}) implies after some
manipulation 
\begin{eqnarray*}
\lefteqn{\hspace*{-1cm}\frac{1}{2k+1}\sum_{t_i\in
    I_c}f(t_i)-\sigma\frac{\sqrt{3\log
      n}+2.72}{\sqrt{2k+1}}}\hspace{4cm}\nonumber\\ 
&\ge&
f(1/2)-\frac{k^2}{2n^2}c_2 -\sigma\frac{\sqrt{3\log
    n}+2.72}{\sqrt{2k}}
\end{eqnarray*}
and on minimizing the right hand side of the inequality with respect to
$k$ we obtain
\begin{eqnarray}
\lefteqn{\hspace*{-2cm}\frac{1}{\vert I_c\vert}\sum_{t_i\in I_c}f(t_i)
  -\sigma\frac{\sqrt{3 \log n}+2.72}{\sqrt{\vert I_c\vert}}}\nonumber\\
&\ge& f(1/2)-1.1c_2^{1/5}\sigma^{4/5}\Big (\sqrt{3\log
  n\,}+2.72\Big)^{4/5}\big/n^{2/5}.\label{requbelow}
\end{eqnarray}
This inequality holds as long as $I_c=[1/2-k_n/n,\, 1/2+k_n/n] \subset
I_0$ with  
\begin{equation} \label{defkn}
k_n = \left\lfloor 0.66c_2^{-2/5}\sigma^{2/5}n^{4/5}\Big(\sqrt{3\log
        n}+2.72\Big)^{2/5}\right\rfloor.  
\end{equation} 
If we put $I_l=[1/2-(\eta+1)k_n/n,\,1/2-\eta k_n/n]$ similar
calculations give 
\begin{eqnarray*}
\lefteqn{\hspace*{-1cm}\frac{1}{2k+1}\sum_{t_i\in I_l}f(t_i)
+\sigma\frac{\sqrt{3\log n}+2.72}{\sqrt{2k+1}}}\hspace{4cm}\\
&\le& f(1/2)-\frac{k^2}{2n^2}c_1 +\sigma\frac{\sqrt{3\log
    n}+2.72}{\sqrt{2k}}
\end{eqnarray*}
and hence
\begin{eqnarray*}
\lefteqn{\hspace*{-1cm}\frac{1}{\vert
I_l\vert}\sum_{t_i\in I_l}f(t_i)+\sigma\frac{\sqrt{3 \log
  n}+´2.72}{\sqrt{\vert I_l\vert}}}\\
&\ge&f(1/2)-\frac{c_2^{1/5}\sigma^{4/5}\Big(\sqrt{3 \log
    n\,}+2.72\Big)^{4/5}}{n^{2/5}}\Big[0.2178\eta^2c_1/c_2-1.23\Big]
\end{eqnarray*}
with the same estimate for $I_r=[1/2+\eta k_n/n,\,1/2+(\eta+1)k_n/n].$
If we put $\eta=3.4 \sqrt{c_2/c_1}$ and  
\begin{equation}
I_n:=\big[1/2-(\eta+1)k_n/n,\,1/2+(\eta+1)k_n/n\big] \subset I_0 
\end{equation}
then all estimates hold. Because of (\ref{defkn}) this will be the case
for $n$ sufficiently large. This implies that (\ref{locmax5}) holds
for sufficiently large $n$ and in consequence any function ${\tilde
  f}_n \in {\mathcal A}_n$ has a local maximum in $I_n.$ 

\subsubsection{Proofs of  (\ref{lbndlocmax0}) and (\ref{ubndlocmax0})}
From (\ref{locmax1}) and (\ref{requbelow}) we have
\[f^*_n(t^*_n)\ge f(1/2)-1.1c_2^{1/5}\sigma^{4/5}\Big (\sqrt{3\log
  n\,}+2.72\Big)^{4/5}\big/n^{2/5}\]
which is the required estimate (\ref{lbndlocmax0}). To prove
(\ref{ubndlocmax0}) we simply note
\[f^*_n(t^*_n)\le f(t^*_n)+\sigma Z(t^*_n)+\sigma\sqrt{3 \log n\,}\le
f(1/2)+\sigma(\sqrt{3 \log n\,}+2.4).\]

\subsubsection{Proof of (\ref{rightadapt}) and (\ref{rate1})}
As $f^*_n \in {\mathcal A}_n$ by definition and $f \in {\mathcal A}_n$ with
probability tending to one we have for the interval $I_{nk}^r=[i/n,(i+k-1)/n]$ 
\[\frac{1}{\sqrt{k}}\sum_{j=0}^{k-1}f_n^*((i+j)/n) \le
\frac{1}{\sqrt{k}}\sum_{j=0}^{k-1}f((i+j)/n)+2\sigma\sqrt{3\log n\,}
\]   
from which it follows that
\[f_n^*(i/n)\le f(i/n)+\frac{k}{n}\Vert
f^{(1)}\Vert_{I_{nk}^r,\infty}+2\sigma\sqrt{\frac{3\log n\,}{k}} \]
which proves (\ref{rightadapt}).
Similarly for the intervals $I_{nk}^l=[(i-k+1)/n,\,i/n]$ we have
\begin{equation} \label{leftadapt}
f(i/n)-f_n^*(i/n)\le \min_{1\le k\le k^{*l}_n}\,
\left\{\frac{k}{n}\Vert f^{(1)}\Vert_{I_{nk}^l,\infty}+
  2\sigma\sqrt{\frac{3\log n\,}{k}}\right\}. 
\end{equation}
We note that (\ref{rightadapt}) and (\ref{leftadapt}) imply that
$f_n^*$ adapts automatically to $f$ to give optimal
rates of convergence. If $f^{(1)}(t) \ne 0$ then it may
be checked that the lengths of the optimal intervals $I_{nk}^{r*}$ and
$I_{nk}^{l*}$ tend to zero and consequently 
\[ \Vert f^{(1)}\Vert_{I_{nk}^{l*},\infty}\approx \vert f^{(1)}(t)\vert
\approx \Vert f^{(1)}\Vert_{I_{nk}^{r*},\infty}. \]
The optimal choice of $k$ is then
\[k^{*l}_n \approx \left(\frac{3\sigma^2 n^2\log n}{\vert
  f^{(1)}(t)\vert^2}\right)^{1/3}\approx k^{*r}_n\]
which gives
\[\lambda(I_{nk}^{l*})\approx \frac{3^{1/3}\sigma^{2/3}}{\vert
  f^{(1)}(t)\vert^{2/3}}\left(\frac{\log n}{n}\right)^{1/3}\approx
\lambda(I_{nk}^{l*})\]
from which (\ref{rate1}) follows.
\newpage
\subsection{Proofs of  Section \ref{concavconvex}}
\subsubsection{Proof of (\ref{lcrconv})}
Then adapting the arguments used
above we have for any differentiable function ${\tilde f}_n \in
{\mathcal A}_n$
\begin{eqnarray*}
\lefteqn{\frac{1}{\sqrt{k}}\sum_{i=1}^{k}\big({\tilde f}_n(1/2+i/n)-
  {\tilde f}_n(1/2-k/n+i/n)\big)}\hspace{2cm}\\ 
&\ge& \frac{1}{\sqrt{k}}\sum_{i=1}^{k}(f(1/2+i/n)-
f(1/2-k/n+i/n))\\
&&\hspace{2cm} -2\sigma(\sqrt{3\log n}+Z(I_{nk}^c)/\sqrt{2})
\end{eqnarray*}
which implies
\begin{equation} \label{centconv}
\max_{t \in I_{nk}^c} f^{*(1)}_n(t)/n \ge \min_{t \in I_{nk}^c}
f^{(1)}(t)/n-\big(2\sigma(\sqrt{3\log n}+
Z(I_{nk}^c)/\sqrt{2})\big)/k^{3/2}.
\end{equation}
Similarly if $I_{nk}^l=[t_l-k/n,t_l+k/n]$ with $t_l+k/n<1/2-k/n$ we
have 
\begin{equation} \label{leftconv}
\min_{t \in I_{nk}^l} f^{*(1)}_n(t)/n \le \max_{t \in I_{nk}^l}
f^{(1)}(t)/n+\big(2\sigma(\sqrt{3\log n}+
Z(I_{nk}^l)/\sqrt{2})\big)/k^{3/2} 
\end{equation}
and for  $I_{nk}^l=[t_r-k/n,t_r+k/n]$ with $t_r-k/n> 1/2+k/n$ we have 
\begin{equation} \label{rightconv}
\min_{t \in I_{nk}^r} f^{*(1)}_n(t)/n \le \max_{t \in I_{nk}^r}
f^{(1)}(t)/n+\big(2\sigma(\sqrt{3\log n}+
Z(I_{nk}^r)/\sqrt{2})\big)/k^{3/2}.
\end{equation}
Again following the arguments given above we may deduce from
(\ref{centconv}), (\ref{leftconv}) and (\ref{rightconv}), that for
sufficiently large $n$ it is possible to choose $I_{nk}^l, I_{nk}^c$
and $I_{nk}^r$ so that (\ref{lcrconv}) holds.

\subsubsection{Proof of (\ref{rightadaptconv})}
We have
\begin{eqnarray*}
\lefteqn{\frac{1}{\sqrt{k}}\sum_{j=1}^k\big(f_n^*(k/n+i/n)-
  f^*_n(i/n)\big)}\\
&\le&\frac{1}{\sqrt{k}}\sum_{j=1}^k\big(f(k/n+i/n)-
  f(i/n)\big)+2\sigma\sqrt{3\log n\,}.
\end{eqnarray*}
and $f^{*(1)}_n$ is non-decreasing on $I^r_{nk}$ we deduce
\[\frac{k^{3/2}}{n}f_n^{*(1)}(t) \le
\frac{1}{\sqrt{k}}\sum_{j=1}^k\big(f(k/n+i/n)- 
  f(i/n)\big)+2\sigma\sqrt{3\log n\,}. \]
A Taylor expansion for $f$ yields
\[f^{*(1)}_n(t) \le f^{(1)}(t)+\frac{k}{n}\Vert
f^{(2)}\Vert_{I^r_{nk},\infty}+2\sigma n \sqrt{\frac{3\log n\,}{k^3}}
\]
from which (\ref{rightadaptconv}) follows.

\subsection{The taut string algorithm of Kovac (2007)}

We suppose that data $y_1,\dots,y_n$ at time points
$t_1<t_2<\cdots<t_n$ are given and first describe how to
calculate the taut string approximation given some tube widths
$\lambda_0,\lambda_1,\dots,\lambda_n$. Subsequently we describe
how to determine these tube widths using a multiresolution
criterion. Lower and upper bounds of a tube on $[0,n]$ are constructed
by linear interpolation of the points $(i,Y_i-\lambda_i),i=0,\dots n$
and $(i,Y_i+\lambda_i),i=0,\dots,n$ respectively where $Y_0=0$ and
$Y_k=Y_{k-1}+y_k$ for $k=1,\dots,n.$ We consider a string $\tilde F_n$
forced to lie in this tube which passes through the points $(0,0)$ and
$(1,Y_n)$ and is pulled tight. An explicit algorithm for doing this
with computational complexity $O(n)$ is described in the
Appendix of Davies and Kovac (2001). The taut string $\tilde F_n$
is linear on each interval $[i-1,i]$ and its derivative
$\tilde f_i=\tilde F_n(i)-\tilde F_n(i-1)$ is used as
an approximation for the data at $t_i$.

Our initial tube widths are $\lambda_0=\lambda_n=0$
and $\lambda_1=\lambda_2=\cdots=\lambda_n=\max(Y_0,\dots,Y_n)-\min(Y_0,
\dots,Y_n)$. We consider the dyadic index set family
\[
{\cal I}=\bigcup_{j,k\in{\Bbb N}_0}\{\{2^jk+1,\dots,2^j(k+1)\}\cap\{1,
\dots,n\}\}\setminus\emptyset
\]
which consists of at most $n$ subsets of $\{1,\dots,n\}$. Given some
taut string approximation $\tilde f_1,\dots,\tilde f_n$ using tube widths
$\lambda_0,\dots,\lambda_n$ we check whether
\begin{equation}\label{mrsum}
\vert\sum_{i\in I} y_i-f_i\vert < \sqrt{3\log(n)}
\end{equation}
is satisfied for each $I\in\cal I$. If this is not the case we generate
new tube widths
$\tilde\lambda_0,\tilde\lambda_1,\dots,\tilde\lambda_n$
by setting $\tilde\lambda_0=\tilde\lambda_n=0$ and for $i=1,\dots,n-1$
\[
\tilde\lambda_i=\begin{cases}
\lambda_i,&\text{if (\ref{mrsum}) is satisfied for all $I\in\cal I$ with
$i\in I$ or $i+1 I$}\\
\lambda_i/2,&\text{otherwise.}\\
\end{cases}
\]
Then we calculate the taut string approximation corresponding to these
new tube widths, check (\ref{mrsum}), possibly determine yet another
set of tube widths and repeat this process until eventually
(\ref{mrsum}) is satisfied for the all $I\in\cal I$.
